\documentclass[12pt,leqno]{amsart}
\usepackage{amsmath,amssymb,amsfonts, amscd,   
pifont, amsthm,  amsxtra, latexsym}
\usepackage[all]{xy}
\usepackage{here}
\usepackage[dvips]{graphicx}

\usepackage{psfrag}
\usepackage{amsmath}
\usepackage{amsthm}
\usepackage{amsfonts}
\usepackage{amssymb}
\theoremstyle{plain}
\newtheorem{thm}{Theorem}[section]

\newtheorem{cor}[thm]{Corollary}
\newtheorem{rem}[thm]{Remark}

\newtheorem{fact}[thm]{Fact}

\newtheorem{remark}[thm]{Remark}

\theoremstyle{definition}
\newtheorem{dfn}[thm]{Definition}

\newtheorem{example}[thm]{Example}

\newcommand{\R}{\mathbb R}
\newcommand{\Z}{\mathbb Z}

\newcommand{\To}{\longrightarrow}

\newcommand{\indhg}{\operatorname{ind}_{\mathfrak h}^{\mathfrak g}}
\newcommand{\indhgprime}{\operatorname{ind}_{\mathfrak h'}^{\mathfrak g'}}
\newcommand{\indpg}{\operatorname{ind}_{\mathfrak p}^{\mathfrak g}}
\newcommand{\indbg}{\operatorname{ind}_{\mathfrak b}^{\mathfrak g}}

\newcommand{\N}{\mathbb N}  %@@
\newcommand{\C}{\mathbb C}  %@@
\makeatletter
\@addtoreset{equation}{section}
\def\@seccntformat#1{\csname the#1\endcsname.\quad}
\def\@seccntformat#1{\csname the#1\endcsname%
\expandafter\ifx\csname#1\endcsname\subsection.\fi\quad}
\makeatother
\newlength{\barlength}
\newcommand{\minusfill}{$\mathsurround=0pt\mathord- \mkern-6mu
   \cleaders\hbox{$\mkern-3mu \mathord- \mkern-3mu$}\hfill
     \mkern-6mu \mathord-$}
\newcommand{\yokobo}{\hbox to 1.2em{\minusfill}}
\newcommand{\equalfill}{$\mathsurround=0pt\mathord= \mkern-6mu
    \cleaders\hbox{$\mkern-3mu \mathord= \mkern-3mu$}\hfill
       \mkern-6mu \mathord=$}
\newcommand{\Longlongrightarrow}
       {\hbox to 2em{\equalfill$\mkern-3mu\Rightarrow$}}
\newcommand{\Longlongleftarrow}
       {\hbox to 2em{$\Leftarrow\mkern-3mu$\equalfill}}
\newcommand{\tsume}{\kern-.35em}
\newlength{\circlength}
\divide\circlength by 2
\advance\circlength by 0.1em

\title{Rankin--Cohen brackets in Representation Theory}
\author[Toshiyuki Kobayashi]{Toshiyuki Kobayashi ${}^{1,2}$}\thanks{${}^{1}$ Graduate School of Mathematical Sciences,
The University of Tokyo, 3-8-1 Komaba, Meguro, Tokyo, 153-8914, Japan}
\thanks{${}^{2}$ French-Japanese Laboratory
 of Mathematics and its Interactions, IRL2025 CNRS and The University of Tokyo, Tokyo, Japan}
\author[Michael Pevzner]{Michael Pevzner ${}^{2,3}$}\thanks{${}^{3}$ 
Universit{\'e} de Reims-Champagne-Ardenne, 
CNRS UMR 9008 LMR, F-51687, Reims, France}

\begin{document}
\maketitle

\begin{abstract}
The Rankin--Cohen brackets provide a basic example of “non-elementary” differential symmetry breaking operators. They can be interpreted as bi-differential operators remarkable for reflecting 
the structure of fusion rules for
holomorphic discrete series representations of the Lie group $SL(2,\mathbb R)$ and are intimately connected to classical special polynomials.   

In this introductory article, we explore the combinatorial structure of these operators and discuss a general framework for constructing their higher-dimensional analogues
from the representa\-tion-theoretic perspective on branching problems. The exposition is based on lectures delivered by the authors during the thematic 
semester ``Representation Theory and Noncommutative Geometry", held in Spring 2025 at the Henri Poincar\'e Institute in Paris. 
\end{abstract}


\medskip
\noindent
\textit{2020 MSC}:
Primary
22E45, % (1973-now) Representations of Lie and linear algebraic groups over real fields: analytic methods {For the purely algebraic theory, see 20G05}
47B38%  (1980-now) Linear operators on function spaces (general),
;
Secondary
11F11, %  (1980-now) Holomorphic modular forms of integral weight
%32A27, % (1980-now) Residues for several complex variables [See also 32C30]
%30H10, %  (2010-now) Hardy spaces [See also 42B30, 46E30]
33C45, %  (1991-now) Orthogonal polynomials and functions of
%hypergeometric type (Jacobi, Laguerre, Hermite, Askey scheme, etc.)
%{For general orthogonal polynomials and functions, see also 42C05}
33C80, %  (1991-now) Connections of hypergeometric functions with groups and algebras, and related topics
43A85. %  (1973-now) Harmonic analysis on homogeneous spaces
%22E45, 22E46, 32M15, 33C45, 33C80, 43A85.
%22E46 (1980-now) Semisimple Lie groups and their representations
%32M15  (1973-now) Hermitian symmetric spaces, bounded symmetric domains, Jordan algebras (complex-analytic aspects) [See also 22E10, 22E40, 53C35, 57T15]

% \noindent{ Declarations of interest: none.}

%\tableofcontents

\section{Introduction}
\label{sec:Intro}

The Rankin--Cohen brackets made their first  implicit appearance in the late 1880s work of P.~Gordan \cite{Gordan} and his student S.~Gundelfinger \cite{Gundel} on invariant theory. Since then, they have been rediscovered on several occasions in various contexts of mathematics, a fact that clearly attests to their enduring relevance. 

Renewed interest in this family of operators largely stems from the role they play in number theory. In a landmark paper \cite{xCo75}, H.~Cohen explored various properties of what were later called
\emph{ the Rankin--Cohen brackets}
 by D.~Zagier \cite{Zagier94},
including an ingenious proof of their covariance, and applied them to the study of special values of $L$-functions.

We therefore introduce these operators through modular forms. 

\subsection{Holomorphic modular forms}
Consider the complex upper half-plane
$$
\mathbb H:=\{z=x+iy~\vert~x\in\mathbb R, y>0\}
$$
endowed with the Poincar\'e metric 
$$
(ds)^2=y^{-2}(dx^2+dy^2)=(\mathrm{Im}(z))^{-2}\vert dz\vert^2
$$ 
that makes $\mathbb H$ into a hyperbolic plane; that is, a real two-dimensional Riemannian manifold of constant negative curvature $-1$.
Moreover, the Poincar\'e metric defines a Hermitian metric compatible with the complex structure. 

The Lie group $SL(2,\mathbb R)$ acts on $\mathbb H$ by fractional linear transformations
\begin{equation}\label{eqn:action}
g:z\mapsto g.z=\frac{az+b}{cz+d},\quad z\in\mathbb H,\, g=\begin{pmatrix}
a&b\\c&d
\end{pmatrix}\in SL(2,\mathbb R).
\end{equation}
%This action is not faithful as $-\mathrm{Id}$ matrix acts trivially. In order to ovoid this one may consider the projective special linear group $PSL(2,\mathbb R):=SL(2,\mathbb R)/\{\pm \mathrm{Id}\}$ instead.
The action of $SL(2,\mathbb R)$ on $\mathbb H$ is transitive; in fact, already the action of its subgroup $P$ of upper triangular matrices is transitive. Indeed,
$$
P\ni\begin{pmatrix} y^{\frac12}&xy^{-\frac12}\\0&y^{-\frac12}\end{pmatrix}
: i\mapsto x+iy.
$$
Moreover, the stabilizer of the point $i\in\mathbb H$ is the subgroup $SO(2)$, which is a maximal compact subgroup of $SL(2,\mathbb R)$. Thus $\mathbb H$ can be viewed as a Hermitian symmetric space isomorphic to $SL(2,\mathbb R)/SO(2)$.
 
  The geometric action \eqref{eqn:action} naturally extends to a family of actions on spaces of functions on $\mathbb H$, parametrized by \emph{weight} $\lambda$.
 
 Consider a particular class of holomorphic functions on $\mathbb H$ that
  satisfy, for a fixed $\lambda\in\mathbb N$, the following functional equation:
\begin{equation}\label{eqn:modular}
 f\left(\frac{az+b}{cz+d}\right)= (cz+d)^\lambda f(z),\quad
\begin{pmatrix}
a&b\\c&d
\end{pmatrix}\in SL(2,\mathbb Z).
 \end{equation}
 
 Since the group $\Gamma=SL(2,\mathbb Z)$ contains the matrix 
 $\begin{pmatrix}
1&1\\ 0&1
\end{pmatrix},
$
any function $f$ satisfying \eqref{eqn:modular} is necessarily periodic, $f(z+1)=f(z)$, and thus admits a Fourier series expansion:
$$
f(z)=\sum_{n\in\mathbb Z}c_n(f)e^{2\pi inz}.
$$
 One says that such a function $f$ is holomorphic at infinity if $c_n(f)=0$ for all $n<0$. If, furthermore, $c_0(f)=0$, then $f$ is said to \emph{vanish at infinity} or to \emph{be cuspidal at infinity}.

\begin{dfn}
Let $\lambda$ be a non-negative even integer. A modular form of weight $\lambda$ with respect to the group $\Gamma=SL(2,\mathbb Z)$ is a holomorphic function
on $\mathbb H$, which satisfies \eqref{eqn:modular} and is holomorphic at infinity.  
\end{dfn}
The vector space of such functions is usually denoted by $M_\lambda(SL(2,\mathbb Z))$ or simply by $M_\lambda(\Gamma )$ if there is no ambiguity. 

 Modular forms exist and form a finite-dimensional vector space
(see, for instance, \cite{Bump} for an elegant introduction to the theory of modular forms). 
For instance, Ramanujan's delta function
$$
\Delta(z)= e^{2\pi iz}\prod_{n=1}^\infty\left(1-e^{2\pi inz}\right)^{24}
$$
is a modular form of weight 12. Moreover, it is cuspidal at infinity.

Given $f\in M_\lambda(\Gamma)$ of weight $\lambda$, it is straightforward to obtain a modular form of weight $2\lambda$ by taking $f^2\in M_{2\lambda}(\Gamma)$, whereas constructing a modular form of
another weight is less immediate. However, the following holds:
\begin{equation}\label{eqn:RC2q}
\lambda f\frac{\partial^2 f}{\partial z^2}-(\lambda+1)\left(\frac{\partial f}{\partial z}\right)^2\in M_{2\lambda+4}(\Gamma).
\end{equation}
It is natural to ask whether there are further formulas like \eqref{eqn:RC2q} to produce modular forms of other higher weights. The answer is affirmative: there exists a countable family of such operators. 

To see this, let us polarize the non-linear operators into bilinear operators. 

We begin with an observation that the holomorphic derivative $\frac{\partial f}{\partial z}$ of a modular form $f$ of weight $\lambda$ is no longer modular, as the second term below constitutes an obstruction :
\begin{eqnarray*}
\frac{\partial f}{\partial z}(z)&=&\frac{\partial}{\partial z}\left((cz+d)^{-\lambda}f\left(\frac{az+b}{cz+d}\right)\right)\\
&=&(cz+d)^{-\lambda-2}\frac{\partial f}{\partial z}\left(\frac{az+b}{cz+d}\right)
-\lambda c(cz+d)^{-\lambda-1}f\left(\frac{az+b}{cz+d}\right).
\end{eqnarray*}
%is no longer modular, since the first term $-\lambda c(cz+d)^{-\lambda-1}f\left(\frac{az+b}{cz+d}\right)$ on the right-hand side fails to satisfy the modularity condition \eqref{eqn:modular}.

However, an appropriate combination of two such holomorphic derivatives of modular forms can cancel the extra terms and produce a \emph{bi-differential operator} that preserves
modularity. For example, let $f_1$ and $f_2$ be two modular forms of weight $\lambda'$ and $\lambda''$, respectively. Then, using the formula above, one observes that the following combination of first-order derivatives:
$$
F(z):=\lambda'f_1(z)\frac{\partial f_2}{\partial z}(z)-\lambda''\frac{\partial f_1}{\partial z}(z)f_2(z)\\
$$
is a modular form of weight $\lambda= \lambda'+\lambda''+2$.
Similarly, the following combination of first- and second-order derivatives:

\begin{equation*}
G(z):=\frac12\lambda''(\lambda''+1)f_1 \frac{\partial^2 f_2}{\partial z^2}
-(\lambda'+1)(\lambda''+1)\frac{\partial f_1}{\partial z}\frac{\partial f_2}{\partial z}+
\frac12\lambda'(\lambda'+1) \frac{\partial^2 f_1}{\partial z^2}f_2
\end{equation*}
is a modular form of weight $\lambda=\lambda'+\lambda''+4$.
 When $\lambda'=\lambda''$ and $f_1=f_2$, the function $F$ vanishes, whereas $G$
reduces to $\eqref{eqn:RC2q}$ 
up to a scalar multiple, and can be viewed as its polarization.

As in the previous observation, this construction is neither \emph{ad hoc} nor isolated.
Rather, it belongs to an infinite family of such bi-differential ``combinations" which
preserve the modularity property \eqref{eqn:modular}. 

\vskip7pt
%Notice that this operator reduces to \eqref{eqn:RC2q} when restricted to the diagonal, that is, when  $f=f_1=f_2$. 

Let us now give a general formula for such transformations.

\begin{dfn}\label{dfn:rcb}
Let $a\in\mathbb N, \lambda',\lambda''\in \mathbb N$ and set $\lambda''':= \lambda'+\lambda''+2a$. 
The $a$-th Rankin--Cohen bracket
$$\mathcal{RC}_{\lambda',\lambda''}^{\lambda'''}\equiv \mathcal{RC}_{\lambda',\lambda''}^{\lambda'+\lambda''+2a}:\mathcal O(\mathbb H)\times\mathcal O(\mathbb H)\to\mathcal O(\mathbb H)
$$ 
is defined by
\begin{multline}\label{eqn:RCB}
\mathcal{RC}_{\lambda',\lambda''}^{\lambda'''}(f_1,f_2)(z):=\\
\mathrm{Rest}_{ z_1=z_2=z}\circ
\sum_{\ell=0}^{a}(-1)^\ell
\left(
        \begin{matrix}
           \lambda'+a-1\\ \ell 
             \end{matrix}
               \right)
               \left(\begin{matrix}
\lambda''+a-1\\a- \ell \end{matrix}\right)
\frac{\partial^{a-\ell} f_1}{\partial z_1^{a-\ell}}(z_1)\frac{\partial^\ell f_2}{\partial z_2^\ell}(z_2).
\end{multline}
Here, $\mathrm{Rest}_{ z_1=z_2
=z}$ denotes the restriction map that simultaneously sets $z_1$ and $z_2$ equal to $z$.
\end{dfn}

For $\lambda'''=\lambda'+\lambda''+2a$ with $a=0,1$ and $2$, the Rankin--Cohen brackets
$\mathcal{RC}_{\lambda',\lambda''}^{\lambda'''}(f_1,f_2)$ take the following form:
\begin{eqnarray*}
&a=0:& f_1f_2,\\
&a=1:& \lambda''\frac{\partial f_1}{\partial z_1}f_2-\lambda'f_1\frac{\partial f_2}{\partial z_1},\\
&a=2:& \frac12\lambda''(\lambda''+1)f_1 \frac{\partial^2 f_2}{\partial z^2}
-(\lambda'+1)(\lambda''+1)\frac{\partial f_1}{\partial z}\frac{\partial f_2}{\partial z}+
\frac12\lambda'(\lambda'+1) \frac{\partial^2 f_1}{\partial z^2}f_2.
\end{eqnarray*}
In particular, the cases $a=0,1,2$ recover the examples discussed earlier.

The Rankin--Cohen brackets satisfy the following fundamental property.
\begin{fact}[Cohen \cite{xCo75}]\label{fact:cova}
Let $a\in\mathbb N, \lambda',\lambda''\in \mathbb N$ and set $\lambda''':= \lambda'+\lambda''+2a$.
If $f_1\in M_{\lambda'}(\Gamma )$ and $f_2\in M_{\lambda''}(\Gamma )$, then $\mathcal{RC}_{\lambda',\lambda''}^{\lambda'''}(f_1,f_2)\in M_{\lambda'''}(\Gamma )$.
\end{fact}

This property reflects the original role of the Rankin--Cohen brackets as an explicit algorithm for constructing modular forms of higher weight from modular forms of lower weight.

Once the explicit expression for the Rankin--Cohen brackets is known, the proof of Fact \ref{fact:cova} reduces to an elementary, though lengthy, direct computation generalizing the preceding observation.
\vskip7pt
Beyond H.~Cohen's original proof \cite{xCo75} of the covariance property in Fact
\ref{fact:cova}, there exist at least seven alternatives proofs based on different ideas, as follows.
\begin{itemize}
\item Recursion relations for the lowest-weight vectors of corresponding representations \cite{Grad, HTan,ABP};
\item Taylor coefficients of Jacobi forms \cite{EiZa};
\item Reproducing kernels of Bergman spaces \cite{vDP, PeZa};
\item Dual pair correspondence \cite{IKO};
\item Source operators \cite{Clerc};
\item Generating operators \cite{KPshort};
\item The F-method \cite{KP16}.
\end{itemize}
In Section \ref{sec:rec}, we will illustrate the first method, which is based on recurrence relations. Within the context of symmetry breaking operators (SBOs), this method is the most elementary. \vskip7pt

On the other hand, to understand the structure of the Rankin--Cohen brackets, we must not only explain the reason why they satisfy the covariance property underlying Fact~\ref{fact:cova}, but also how they can be discovered anew and why, up to a scalar, no other such bi-differential operators exist provided that $\lambda',\lambda''\in\N$. To achieve this, we appeal to the representation theory of the ambient Lie group $SL(2,\mathbb R)$, rather than viewing the covariance as a purely arithmetic feature of the modular group $SL(2,\mathbb Z)$.

The last method, referred to as the F-method, derives the explicit form of the Rankin--Cohen brackets from scratch by revealing that the seemingly complicated alternating sum coefficients arise naturally from the Jacobi polynomials. Moreover, it provides ``exhaustion-type" theorems in the classification theory of SBOs, see Corollary \ref{cor:83}. 

The F-method will be discussed in a more general framework in Section \ref{sec:2}, based on the ``algebraic
Fourier transform of Verma modules". Meanwhile, we introduce some  basic notions that will be used later for this purpose.

\subsection{Holomorphic discrete series representations}
Among irreducible unitary representations of a real simple Lie group $G$, there exists a distinguished family of infinite-dimensional representations, that occurs as subrepresentations of the left regular representation of $G$ on the Hilbert space $L^2(G)$. These representations are called \emph{square-integrable representations}, because they are characterized by the square-integrability of their matrix coefficients. They are also referred to as \emph{discrete series representations} of $G$.
Harish-Chandra initiated a systematic study of these representations and provided both a criterion for their existence in terms of the rank condition, as well as a parametrization of discrete series representations. See Duflo \cite{Duflo} for a survey of the developments of the theory of discrete series during the 1970s.\vskip7pt

The group $SL(2,\mathbb R)$ possesses a non-empty discrete series. In this case, the discrete series exhibits a particularly rich structure: the corresponding $(\mathfrak g,K)$-modules are either lowest or highest-weight modules. 
 These discrete series representations can be realized concretely in Hilbert spaces of holomorphic or anti-holomorphic functions respectively. 
 
 More precisely, for $\lambda\in\mathbb Z$, consider the space 
$$
\mathcal H_\lambda^2(\mathbb H):=\mathcal O(\mathbb H)\cap L^2(\mathbb H, y^{\lambda-2}dxdy),
$$
consisting of holomorphic functions on the upper half-plane $\mathbb H$ that are square integrable with respect to the weighted measure $ y^{\lambda-2}dxdy$. If $\lambda>1$, then $\mathcal H_\lambda^2(\mathbb H)$ is nonzero and is infinite-dimensional. Moreover, since the pointwise evaluation map is a continuous linear functional on this  Hilbert space, it admits a reproducing kernel. The space $\mathcal H_\lambda^2(\mathbb H)$ is referred to as a \emph{weighted Bergman space}.

The fractional linear action \eqref{eqn:action} of the group $SL(2,\mathbb R)$ on $\mathbb H$ induces a unitary representation
$\pi_{\lambda}$ on $\mathcal H_\lambda^2(\mathbb H)$, given by:
\begin{equation}\label{eqn:hds}
(\pi_\lambda(g)f)(z)=(cz+d)^{-\lambda}f\left(\frac{az+b}{cz+d}\right),\,
{\mathrm {for}} \,g^{-1}=\begin{pmatrix}
a&b\\c&d
\end{pmatrix}\in SL(2,\mathbb R).
 \end{equation}

For $\lambda>1$ such representations $\pi_\lambda$ are irreducible and pairwise inequivalent. The fact that these representations can be realized in the kernel of the Cauchy--Riemann operator is a special feature of the group $SL(2,\mathbb R)$;
representations arising in this manner are referred to as \emph{holomorphic discrete series} representations. More generally, this phenomenon persists in the broader framework of Hermitian Lie groups. 

We note that the classical Paley--Wiener theorem provides yet another realization of the \emph{holomorphic discrete series representations}
of $SL(2,\mathbb R)$. Indeed, the Fourier--Laplace transform maps the weighted Bergman space $\mathcal H^2_\lambda(\mathbb H)$ into the Hilbert space $L^2(\R_+, x^{1-\lambda}d x)$, thus
 giving rise to an
$L^2$-\emph{model} of the same representation $\pi_\lambda$. Although explicit formulas for the group action are more complicated than in
\eqref{eqn:hds}, this model turns out to be appropriate for the analysis of Rankin--Cohen brackets and the dual notion of symmetry breaking operators, namely the \emph{ holographic operators} \cite{KP20}, through the F-method.

We also note that the action \eqref{eqn:hds} admits a more conceptual interpretation. Indeed, there exists a countable family of $SL(2,\mathbb R)$-equivariant holomorphic line bundles $\mathcal L_\lambda$ over $\mathbb H$, 
parametrized by $\lambda\in\mathbb Z$, associated with a character of $SO(2)$, which together exhaust all such bundles. The regular representation of the group $SL(2,\mathbb R)$ on the space of holomorphic sections $\mathcal O(\mathbb H,\mathcal L_\lambda)$ can be identified with $\pi_\lambda$ after the trivialization of this line bundle via the Bruhat decomposition. 

The covariance property of the Rankin--Cohen brackets in Definition \ref{dfn:rcb} is formulated in terms of the representations
$\pi_\lambda$ given in \eqref{eqn:hds} as follows.
\begin{fact}
Let $a\in\mathbb N, \lambda',\lambda''\in \mathbb N$ and set $\lambda''':= \lambda'+\lambda''+2a$.
Then for every $f_1,f_2\in \mathcal O(\mathbb H)$ and any $g\in SL(2,\mathbb R)$ one has 
$$
\mathcal{RC}_{\lambda',\lambda''}^{\lambda'''}(\pi_{\lambda'}(g)f_1,\pi_{\lambda''}(g)f_2)(z)=\left(\pi_{\lambda'''}(g)\mathcal{RC}_{\lambda',\lambda''}^{\lambda'''}(f_1,f_2)\right)(z).
$$
\end{fact}
This result extends Fact \ref{fact:cova} to a larger group of transformations---namely, from 
$SL(2,\mathbb Z)$ to $SL(2,\mathbb R)$---and shows that the covariance property is fundamentally geometric rather than arithmetic holding over the entire Lie group $SL(2,\mathbb R)$. 
Later, we allow the
parameter $\lambda$ to be a complex number, and one observes that this statement remains valid for the universal covering group of ${SL}(2,\mathbb R)$.

Moreover, it admits a representation-theoretic interpretation, as the Rankin--Cohen brackets can be regarded as intertwining operators, see for instance \cite{vDP}. We adopt this perspective and demonstrate how it explains the structure of these operators, why they take this specific form, and how to construct, in an exhaustive and systematic way, their analogues for different families of real semisimple Lie groups.

\subsection{Weight Watchers}\label{sec:rec}

One somewhat handcrafted way to construct such intertwining operators relies on the fact that the holomorphic discrete series representations are lowest-weight modules satisfying weight-multiplicity-free property. This property provides an elementary method for ${SL}(2,\mathbb R)$, yet it is not easily applicable for higher-rank groups.
\vskip10pt

Consider the real Lie algebra of $2\times 2$ traceless matrices 
$$\mathfrak{sl}(2,\mathbb R):=\{X\in \mathrm{Mat}(2,\mathbb R):\mathrm{tr}(X)=0\}.$$ 

As a vector space, it is three -dimensional, and we choose the following basis elements
$$
e^+= \begin{pmatrix}
0&1\\0&0
\end{pmatrix},
\quad
h=\begin{pmatrix}
1&0\\0&-1
\end{pmatrix},
\quad
e^-=\begin{pmatrix}
0&0\\1&0
\end{pmatrix},
$$
which are subject to the commutation relations:
$$
[h,e^\pm]=\pm2e^\pm,\quad [e^+,e^-]=h.
$$
The center $\mathcal Z$ of the universal enveloping algebra $U(\mathfrak{sl}_2)$ of the complexification of the Lie algebra $\mathfrak{sl}(2,\mathbb R)$ is generated, as a $\C$-algebra, by the single element:
$$
c=h^2+2(e^+e^-+e^-e^+)=h^2+2h+4e^-e^+.
$$

The infinitesimal action $d\pi_\lambda$ of $\mathfrak{sl}(2,\mathbb R)$ on the space of holomorphic functions $\mathcal O(\mathbb H)$ corresponding to the holomorphic discrete series representation $\pi_\lambda$ of the Lie group $SL(2,\mathbb R)$ is defined by
$$
d\pi_\lambda(X)f (z) := \frac{d}{dt}\big\vert_{t=0}\pi_\lambda(\exp(tX)) f(z),\, X\in \mathfrak{sl}(2,\mathbb R), f\in\mathcal O(\mathbb H).
$$

Then, according to the formula \eqref{eqn:hds} we have
\begin{eqnarray*}
d\pi_\lambda(e^-)f(z) &=& \lambda z f(z)+ z^2 f'(z),\\
d\pi_\lambda(e^+)f(z) &=&-f'(z),\\
d\pi_\lambda(h)f(z) &=&-\lambda f(z)-2zf'(z).
\end{eqnarray*}

Furthermore, consider the particular element 
\begin{equation}\label{eqn:hwv}
{}^\lambda v_\ell\equiv v_\ell:=\frac{(\lambda+\ell-1)!}{(\lambda-1)!}z^{-\lambda-\ell}
\end{equation} 
and
denote by $V_\lambda=\mathrm{Span}_{\C}\langle {}^\lambda v_\ell\rangle_{\ell\in\mathbb N}$ the vector space spanned by these functions. Then,
\begin{eqnarray*}
d\pi_\lambda(e^+)({}^\lambda v_\ell)&=&{}^\lambda v_{\ell+1}, \quad\forall\ell\in\mathbb N,\\
d\pi_\lambda(e^-)({}^\lambda v_\ell)&=&-\ell(\lambda+\ell-1)\,{}^\lambda v_{\ell-1}, \quad\forall\ell\in\mathbb N\setminus\{0\},\\
d\pi_\lambda(e^-)({}^\lambda v_0)&=&0,\\
d\pi_\lambda(h)({}^\lambda v_\ell)&=&(\lambda+2\ell){\,}^\lambda v_{\ell}, \quad\forall\ell\in\mathbb N,\\
d\pi_\lambda(c)({}v)&=&(\lambda^2-2\lambda)v, \quad\forall v\in V_\lambda,
\end{eqnarray*}
where in the last identity, we abuse notation slightly by using the same symbol for the extension of the Lie algebra representation $d\pi_\lambda$ to its universal enveloping algebra.

Thus, the space $V_\lambda$ can be seen as a lowest-weight module for the Lie algebra $\mathfrak{sl}(2,\mathbb R)$. Moreover, the maximal compact subgroup $SO(2)$  of $SL(2,\mathbb R)$ acts on $V_\lambda$ in such a way that every
element of $V_\lambda$ is $SO(2)$-finite in the sense that 
$\dim \mathrm{Span}_{\C}\langle \pi_\lambda(g)v : g \in SO(2)\rangle  < \infty$ for all $v\in V_\lambda$. 
The space $V_\lambda$ is dense in $\mathcal H^2_\lambda(\mathbb H)$, and $V_\lambda$ is the underlying $(\mathfrak{g}, K)$-module for the holomorphic discrete series representation $\pi_\lambda$. 
Since $d\pi_\lambda(e^-)({}^\lambda v_0)=0$, the Harish-Chandra module $V_\lambda$ is a lowest-weight module.

For two given integers $\lambda',\lambda'' >1$, consider the tensor product $V_{\lambda'}\otimes V_{\lambda''}$ of two such modules as an $\mathfrak{sl}(2,\mathbb R)$-module via the diagonal embedding of the Lie algebra $\mathfrak{sl}(2,\mathbb R)$ into $\mathfrak{sl}(2,\mathbb R)\oplus \mathfrak{sl}(2,\mathbb R)$. Then this module is no longer irreducible and decomposes discretely and multiplicity-freely into a direct sum of lowest-weight Harish-Chandra modules:
\begin{equation}\label{eqn:bra-vk}
V_{\lambda'}\otimes V_{\lambda''}=\bigoplus_{a\in\mathbb N} V_{\lambda'+\lambda''+2a}.
\end{equation}
The $\mathfrak{sl}(2,\mathbb R)$-case was proved by Molchanov \cite{Mol79} and Repka \cite {Repka}.  We refer to \cite[Theorem 8.4]{Kmf} for a more conceptual explanation of this phenomenon in a broader framework.

As a refinement of the abstract branching rule \eqref{eqn:bra-vk}, one can explicitly determine the lowest-weight vectors of the irreducible components appearing in the decomposition \eqref{eqn:bra-vk} (see \cite{Grad,HTan}, for example). 

We now carry out this computation. The element $\left({}^{\lambda'}v_\ell\right)\otimes\left({}^{\lambda''}v_m\right)\in V_{\lambda'}\otimes V_{\lambda''}$ is an eigenvector for the diagonal action of $h\in\mathfrak{sl}(2,\mathbb R)$, with eigenvalue $\lambda'+\lambda''+2(\ell+m)$.

Therefore, the eigenspace corresponding to the eigenvalue 
$$
\lambda''':=\lambda'+\lambda''+2a
$$ 
is generated, as a vector space, by elements
of the form ${}^{\lambda'}v_\ell\otimes{}^{\lambda''}v_{a-\ell}$.\vskip7pt

To identify irreducible components in the decomposition \eqref{eqn:bra-vk}, we thus need to find a family of constants $\kappa_\ell\in\mathbb C$ such that  the element
$$
w_{\lambda'\lambda''}^{a}(z_1,z_2):=\sum_{\ell=0}^{a}\kappa_\ell ({}^{\lambda'}v_\ell\otimes{}^{\lambda''}v_{a-\ell})\in V_{\lambda'}\otimes V_{\lambda''}
$$
is annihilated by the diagonal action  of the element $e^-\in\mathfrak{sl}(2,\mathbb R)$.

Thus, in order to determine these coefficients, we describe the kernel of the differential operator 
\begin{eqnarray*}
d\pi_{\lambda'}(e^-)\otimes \mathrm{Id}_{V_{\lambda''}}&+&
 \mathrm{Id}_{V_{\lambda'}}\otimes d\pi_{\lambda''}(e^-)\\
=\left(\lambda'z_1+z_1^2\frac{\partial}{\partial z_1}\right)\otimes \mathrm{Id}_{V_{\lambda''}}&+&
\mathrm{Id}_{V_{\lambda'}}\otimes 
\left(\lambda''z_2+z_2^2\frac{\partial}{\partial z_2}\right),
\end{eqnarray*}
where the variables $z_1, z_2$ refer to the realization of the modules $V_{\lambda'}$ and $V_{\lambda''}$, respectively.

If $w_{\lambda'\lambda''}^{a}$ is annihilated by this operator, the comparison of the coefficients of 
${}^{\lambda'}v_\ell\otimes{}^{\lambda''}v_{a-\ell}$ yields the following equality for every $\ell$:
\begin{equation}\label{eqn:recur}
(\ell+1)(\lambda'+\ell)\kappa_{\ell+1}+(a-\ell)(\lambda''+a-\ell-1)\kappa_\ell=0,
\end{equation}
Interpreting this as a recurrence relation, one has
\begin{equation}\label{eqn:lambda}
\kappa_{\ell}=\frac{a!}{(\lambda'')_a}(-1)^\ell
\left( \begin{matrix}
           \lambda'+a-1\\ \ell 
             \end{matrix}
               \right)
               \left(\begin{matrix}
\lambda''+a-1\\a- \ell \end{matrix}\right)
 \kappa_0,
\end{equation}
where $(x)_n=\frac{\Gamma(x+n)}{\Gamma(x)}$ is the rising Pochhammer symbol and the initial condition $\kappa_0$ can be chosen freely.

%Notice that, equivalently, one may interpret the above condition as a system of linear equations of size $a\times (a+1)$ of very specific form, see \cite{Grad}), that leads, of course, to the same solution.

%Let us choose $\kappa_0$ in such a way that
%\begin{eqnarray}\nonumber
%%:
%w_{\lambda'\lambda''}^{a}(z_1,z_2)&=&\sum_{\ell=0}^a(-1)^\ell
%\begin{pmatrix}
%a\\
%\ell
%\end{pmatrix}
%\frac{(\lambda'+a-1)!(\lambda''+a-1)!}{a!(\lambda'-1)!(\lambda''-1)!}z_1^{-\lambda'-\ell}z_2^{-\lambda''-a+\ell}\\
%&=&\frac{(\lambda'+a-1)!(\lambda''+a-1)!}{a!(\lambda'-1)!(\lambda''-1)!}
%\frac{(z_2^{-1}-z_1^{-1})^a}{z_1^{\lambda'}z_2^{\lambda''}}.
%\label{eqn:phi}
%\end{eqnarray}
The element $w_{\lambda'\lambda''}^{a}$ is, up to a scalar multiple, the lowest-weight vector in $V_{\lambda'''}= V_{\lambda'+\lambda''+2a}$, viewed as an irreducible submodule of $V_{\lambda'}\otimes V_{\lambda''}$.

Since $d\pi_\lambda(e^+)({}^\lambda v_\ell)={}^\lambda v_{\ell+1}$, the element $w_{\lambda'\lambda''}^{a}$ can also be viewed as the image of the tensor product of the two lowest-weight vectors
$\left({}^{\lambda'}v_0\right)\otimes\left({\,}^{\lambda''}v_0\right)$ under the action of the element
$$
\mathfrak{RC}_{\lambda'\lambda''}^{\lambda'''}:= \sum_{\ell=0}^a
 \kappa_\ell
(e^+)^{\ell}\otimes (e^+)^{a-\ell}\in U(\mathfrak{sl}_2)\otimes U(\mathfrak{sl}_2).
$$

Duality Theorem \ref{thm:surject} in Section \ref{sec:2} asserts that finding a singular vector (\emph{i.e.} a lowest-weight vector in the tensor product representation) is equivalent to constructing  a differential symmetry breaking operator.
In the present setting,
the corresponding bi-differential operator is obtained by formally replacing
 $e^+$ with $dR(e^+)=\frac{\partial}{\partial z}$ (see Example \ref{exa:1}).
 Consequently, the element
  $\mathfrak{RC}_{\lambda'\lambda''}^{\lambda'''}$
  in the enveloping algebra induces a bi-differential operator
   on $\mathcal O(\mathbb H\times\mathbb H)$ given by
\begin{eqnarray*}
&&\sum_{\ell=0}^a
 \kappa_\ell
 \left(\frac{\partial}{\partial z_1}\right)^\ell
 \left(\frac{\partial}{\partial z_2}\right)^{a-\ell}\\
&=&\kappa_0\frac{(-1)^aa!}{(\lambda'')_a}
\sum_{\ell=0}^a
(-1)^{a-\ell}
\left(
        \begin{matrix}
           \lambda'+a-1\\ \ell 
             \end{matrix}
               \right)
               \left(\begin{matrix}
\lambda''+a-1\\a- \ell \end{matrix}\right)
\frac{\partial^a}{\partial z_1^{\ell}\partial z_2^{a-\ell}}
\end{eqnarray*}
by \eqref{eqn:lambda}. Restricting this operator to $z_1=z_2(=z)$,
  and choosing the normalizing constant 
  $$
  \kappa_0=\frac{(-1)^a(\lambda'')_a}{a!},
  $$ 
 we recover precisely the $a$-th Rankin--Cohen bracket $\mathcal{RC}_{\lambda'\lambda''}^{\lambda'''}$. 
Therefore, by construction, we verify that $\mathcal{RC}_{\lambda'\lambda''}^{\lambda'''}\in\mathrm{Hom}_{\mathfrak{sl}(2,\mathbb R)}(V_{\lambda'}\otimes
V_{\lambda''}, V_{\lambda'''})$.

\subsection{Connection to orthogonal polynomials}

Far from being an isolated instance, this example fits into a broader framework, developed in
Section \ref{sec:2}, which clarifies the combinatorial structure of such homomorphisms and opens the door to far-reaching generalizations.

At first glance, one might expect that taking normal derivatives with respect to equivariant embeddings of homogeneous spaces would yield intertwining operators (symmetry breaking operators) valued in representations of the corresponding subgroups acting on functions defined on these submanifolds. Although this occurs only under restrictive conditions, one encounters phenomena in which
 intertwining operators appearing in branching rules are realized by holomorphic differential operators that cannot be expressed as normal derivatives (see \cite{KP16}, Section 5). 

The specific geometric structure of the upper half-plane and of its higher-rank analogues constrains
the class of intertwining operators under consideration to be local operators, as formulated by Localness Theorem \ref{thm:C} below.

Rankin--Cohen brackets \eqref{eqn:RCB} are holomorphic differential operators with constant coefficients. As seen in \eqref{eqn:lambda}, the coefficients obtained from the recurrence relation, 
coincides, up to normalization, exactly with those of certain orthogonal polynomials. 
Indeed, using the notation introduced in \cite{{KP16}},
the Rankin--Cohen brackets  may be expressed as
\begin{equation}\label{eqn:rcJacob}
\mathcal{RC}_{\lambda'\lambda''}^{\lambda''}=(-1)^a
\mathrm{Rest}_{z=z_1=z_2}
\circ P_a^{\lambda'-1,1-\lambda'''}\left(\frac{\partial}{\partial z_1},\frac{\partial}{\partial z_2}\right),
\end{equation}
where $ P_a^{\alpha,\beta}(x,y)$ is a homogeneous polynomial in two variables of degree $a$,
defined by
\begin{equation}\label{eqn:Jacobi_two}
P_a^{\alpha,\beta}(x,y):=y^a P_a^{\alpha,\beta}\left(2\frac{x}{y}+1\right),
\end{equation}
and obtained by ``inflation" (see \eqref{eqn:inflJacobi}, Section 2) from the classical Jacobi polynomial, which is given by the
Rodrigues formula
\begin{equation}\label{eqn:rodrigues}
P_a^{\alpha,\beta}(t) =\frac{(-1)^a}{2^aa!}(1-t)^{-\alpha}(1+t)^{-\beta}\left(\frac{d}{dt}\right)^a\left(
(1-t)^{a+\alpha}(1+t)^{a+\beta}\right).
\end{equation}

This precise correspondence between the Rankin--Cohen brackets and the Jacobi polynomials is not accidental, but can be explained systematically by the F-method discussed in Section \ref{sec:2}. Within this framework, an important aspect is
that the Jacobi polynomial $P_a^{\alpha,\beta}(t)$ arises as a polynomial solution of the Gauss hypergeometric differential equation (but not vice versa, see Theorem \ref{thm:Jacobizero}).
To make this explicit,
expressing the Jacobi polynomials in terms of the Gauss hypergeometric function
$$
P_a^{\alpha,\beta}(t)=\frac{(\alpha+1)_a}{a!}{}_2F_1\left(-a;\alpha+\beta+a+1,\alpha+1;\frac{1-t}2\right),
$$
one finds that they satisfy the following second-order differential equation, known as
the \emph{Jacobi differential equation}:
\begin{equation}\label{eqn:JacobiDE}
\left((1-t^2)\frac{d^2}{dt^2}+(\beta-\alpha-(\alpha+\beta+2)t)\frac d{dt}+a(a+\alpha+\beta+1)\right)y=0.
\end{equation}

The \emph{Rodrigues formula} \eqref{eqn:rodrigues} provides a recursive structure for the entire family of Jacobi polynomials and ensures that they form a complete orthonormal system of the Hilbert space $L^2([-1,1],(1-t)^{\alpha}(1+t)^{\beta}dt)$. These three key features of Jacobi polynomials---differential equation, recursion, and orthogonality---acquire a profound significance in the analysis of branching rules.

While the role of orthogonality in the analysis of symmetry breaking operators, or their dual notion of holographic operators was first addressed in \cite{KP20}, and will not be pursued here, we shall show in the next section that the differential equation aspect arises naturally via the ``algebraic Fourier transform of Verma modules".

Before delving further into the representation-theoretic perspective, let us note that beyond their role in the theory of modular forms, the Rankin--Cohen brackets also appear in a variety of other settings, including equivariant deformation quantization, as well as non-commuta\-tive and conformal geometry \cite{CMZ,CM,Juhl,UU96}.

\section{Symmetry breaking Paradigm}\label{sec:2}

\emph{Symmetry breaking in branching problems} offers a natural and conceptually transparent framework for understanding Rankin--Cohen bra\-ckets and their higher-dimensional  analogues. 
In this section, we outline this perspective and explain why it is both natural and effective.
Our exposition follows the approach developed in \cite{KP16a,KP16}. 

Let
$p:Y\To X$ be a smooth map between two manifolds $Y$ and $X$.  {Let $\mathcal W\To Y$ and $\mathcal V\To X$ be vector bundles over $Y$ and $X$, respectively. 
A central object in this framework is a linear operator relating sections over $X$ to those over $Y$. 
Following \cite{KP16a}, we say that a continuous linear map $T:C^\infty(X,\mathcal V)\to C^\infty(Y,\mathcal W)$ is a \emph{differential
operator between two manifolds} if it satisfies the natural support condition
\begin{equation}
\label{eqn:suppT}
p\left(\operatorname{Supp} T f\right)
 \subset
 {\operatorname{Supp}f}
 \end{equation}
 for all smooth sections
$ f\in C^\infty (X,\mathcal V)$.

The support condition \eqref{eqn:suppT} implies
 the \emph{localness} of the operator $T$ along the map $p$ 
 in the following sense: for any open subset $U$ of $X$, the operator
$T$ induces a continuous linear map:
$$
T_U: C^\infty(U,\mathcal V\vert_U)\longrightarrow C^\infty\left(p^{-1}(U),\mathcal W\vert_{p^{-1}(U)}\right).
$$

If $X=Y$
 and $p$ is the identity map,
 then the condition \eqref{eqn:suppT} is equivalent 
 to $T$ being a differential operator
 in the usual sense,
 by Peetre's theorem \cite{Pee59}. 
 
 In the general case, any operator $T$ satisfying \eqref{eqn:suppT} can be expressed as a finite linear
 combination of operators of the form $F_Y\circ p^*\circ D_X$, where $F_Y$ denotes
 multiplication by a smooth function on $Y$, and $D_X$ is a differential operator on $X$.
 A characterization of such operators in terms of their Schwartz distribution kernels was given in \cite[Lemma 2.3]{KP16a}.
 
 We denote by $\operatorname{Diff}(\mathcal V_X,\mathcal W_Y)$ the vector space
 of differential operators
from $C^\infty(X,\mathcal V)$ to $C^\infty(Y,\mathcal W)$.

 This framework plays a fundamental role in the formulation of differential symmetry breaking operators. With this in mind, we now incorporate symmetries into the setting
 by considering a pair of Lie groups $G'\subset G$.
We assume that $G'$ acts on $Y$, that $G$ acts on $X$, and
that the map $p:Y\To X$ is $G'$-equivariant. Moreover, we assume that vector bundles $\mathcal W\To Y$ and $\mathcal V\To X$ are equivariant with respect to the actions of $G'$ and $G$, respectively.

 The Lie group $G'$ acts on the vector bundle $\mathcal V\to X$ by restriction of the action of $G$.
 This action induces natural actions of $G'$ on the Fr{\'e}chet spaces
 $C^\infty(X,\mathcal{V})$
and $C^\infty(Y,\mathcal{W})$ by translations.
Denote by $\operatorname{Hom}_{G'}(C^\infty(X,\mathcal V),C^\infty(Y,\mathcal W))$
 the space of continuous $G'$-homomor\-phisms and introduce the main object of our study: 
\begin{eqnarray*}
\operatorname{Diff}_{G'}(\mathcal V_X,\mathcal W_Y):=
\operatorname{Diff}(\mathcal V_X,\mathcal W_Y)\cap
\operatorname{Hom}_{G'}(C^\infty(X,\mathcal V),C^\infty(Y,\mathcal W)).
\end{eqnarray*}

We refer to elements of $\operatorname{Hom}_{G'}(C^\infty(X,\mathcal V), C^\infty (Y,\mathcal W))$ as \emph{symmetry breaking operators}, and to elements of $\operatorname{Diff}_{G'}(\mathcal V_X,\mathcal W_Y)$ as
 \emph{differential symmetry breaking operators}.

To illustrate these notions, consider the following example. If $X$ and $Y$ are both Euclidean vector spaces
 with an injective linear map $p: Y \hookrightarrow X$ and if
$G'$ contains the subgroup of all translations of $Y$, then
 $\operatorname{Diff}_{G'}(\mathcal V_X,\mathcal W_Y)$ is a subspace of the space of
differential operators with constant coefficients.

If the vector bundle $\mathcal W\To Y$ is isomorphic to the pullback $p^*\mathcal V$ of $\mathcal V$, then the restriction map $f\mapsto f\vert_Y$ is a $G'$-equivariant differential operator of order zero. 

An analogous notion can be defined
 in the holomorphic setting.  
Let ${\mathcal {V}} \to X$ and ${\mathcal {W}} \to Y$
 be two holomorphic vector bundles with a holomorphic map $p:Y\to X$ between the complex
manifolds $X$ and $Y$.  
We say that a differential operator $T: C^\infty(X,\mathcal V)\to C^\infty(Y,\mathcal W)$ is {\it{holomorphic}} if
$$
T_U(\mathcal O(U,\mathcal V\vert_U))
 \subset\mathcal O(p^{-1}(U),\mathcal W\vert_{p^{-1}(U)})
$$
for any open subset $U$ of $X$. We denote by $\operatorname{Diff}^{\mathrm{hol}}(\mathcal V_X,\mathcal W_Y)$ the vector space of holomorphic differential operators.
When a Lie group $G'$ acts biholomorphically on the two holomorphic vector bundles $\mathcal V\to X$ and $\mathcal W\to Y$, we set
$$
\operatorname{Diff}_{G'}^{\mathrm{hol}}(\mathcal V_X,\mathcal W_Y):=
\operatorname{Diff}^{\mathrm{hol}}(\mathcal V_X,\mathcal W_Y)
\cap
\operatorname{Hom}_{G'}(C^\infty(X,\mathcal V), C^\infty (Y,\mathcal W)).
$$

Symmetry breaking operators arising from concrete geometric models provide both a guiding framework and powerful tools for
the detailed study of branching problems.

A new method, the F-method, based on the \emph{algebraic Fourier transform of generalized Verma modules}, was developed in \cite{KP16a} to characterize the covariance condition in $\operatorname{Diff}_{G'}(\mathcal V_X,\mathcal W_Y)$ in terms of systems of higher-order differential equations.
It provides both an explicit construction and an effective tool for proving the exhaustion
 of differential symmetry breaking operators. This approach is based on three key results, which we now describe.

 \subsection{The big three}

Suppose that the groups $G$ and $G'$ act transitively on $X$ and $Y$, respectively, so that we may write $X=G/H$ and $Y=G'/H'$.

The first result, referred to as the \emph{Duality Theorem} (\cite[Theorem 2.9]{KP16a}), provides an algebraic description of symmetry breaking operators in terms of generalized Verma modules:

\begin{thm}[Duality Theorem]\label{thm:surject}
Let $H' \subset H$ be (possibly disconnected) closed subgroups
 of a Lie group $G$ with Lie algebras
$\mathfrak{h}' \subset \mathfrak{h}$,
respectively.
Suppose  that 
$V$ and $W$ are finite-dimensional  representations of 
$H$ and $H'$, respectively, and denote by $V^\vee$ and $W^\vee$ corresponding contragredient representations.

Let $G'$ be any subgroup of $G$ containing $H'$, and $\mathcal V_X:=G\times_{H} V$ and 
 $\mathcal W_Y:=G'\times_{H'} W$ be the corresponding
homogeneous
vector bundles. 
Let $\indhg(V^\vee)$ be the $\mathfrak g$-module defined by
$
\indhg(V^\vee):=U(\mathfrak g)\otimes_{U(\mathfrak h)} V^\vee.
$ Notice that this is a generalized Verma module if $\mathfrak h$ is a parabolic subalgebra.
Define $dR\colon\mathfrak g\To\mathrm{End}(C^\infty(G,V))$ by
 $$
   \left(dR(X)F\right)(g)=\frac{d}{dt}\big\vert_{t=0}F(ge^{tX}).
 $$
 Then there is a natural linear isomorphism:
\begin{equation}
\label{eqn:DXYH}
D_{X\to Y}:
\operatorname{Hom}_{H'}(W^\vee,\indhg(V^\vee))
\stackrel{\sim}{\longrightarrow}
\operatorname{Diff}_{G'}\left( \mathcal V_X,
\mathcal W_Y\right),
\end{equation}
or equivalently,
$$
D_{X\to Y}:
\operatorname{Hom}_{(\mathfrak g',H')}(\indhgprime(W^\vee),\indhg(V^\vee))
\stackrel{\sim}{\longrightarrow}
\operatorname{Diff}_{G'}\left( \mathcal V_X,
\mathcal W_Y\right).
$$
For $\varphi \in \operatorname{Hom}_{H'}(W^\vee,\indhg(V^\vee))$
 and $F \in C^{\infty}(X, {\mathcal{V}}) \simeq C^{\infty}(G,V)^H$, 
 the corresponding section
 $D_{X \to Y}(\varphi)F \in C^{\infty}(Y, {\mathcal{W}}) \simeq C^{\infty}(G',W)^{H'}$
 is defined by the following formula:
\begin{equation}
\label{eqn:Dxy}
\langle 
D_{X \to Y}(\varphi)F, 
w^{\vee}
\rangle
=
\sum_j
\langle 
d R(u_j)F
, 
v_j^{\vee}
\rangle
|_{G'}
\quad
\text{for $w^{\vee} \in W^{\vee}$}, 
\end{equation}
where $u_j \in U({\mathfrak {g}})$ and $v_j^{\vee} \in V^{\vee}$ are chosen in such a way that
 $\varphi(w^{\vee})=\sum_j u_j v_j^{\vee}$
 $\in \indhg(V^\vee)$. 
 \end{thm}

\begin{example}\label{exa:1}
In the setting where $V=W=\C$, the homomorphism $\varphi\in \operatorname{Hom}_{H'}(W^\vee,\indhg(V^\vee))$
is represented by $u\in U(\mathfrak g)$ such that $\varphi(1)=u\cdot 1$. Consequently, the formula
\eqref{eqn:Dxy} reduces to
$$
D_{X \to Y}(\varphi)F= dR(u)F\big\vert_{G'}. 
$$
\end{example}

When $H'$ is connected,
the left-hand side of \eqref{eqn:DXYH} can be expressed in terms of Lie algebras.  
\begin{cor}
\label{cor:2.9}
Suppose we are in the setting  of
 Theorem \ref{thm:surject} and
 that $H'$ is connected.  
Then the morphism \eqref{eqn:DXYH} is bijective :
\begin{eqnarray}\label{eqn:isoo}
D_{X\to Y}:\operatorname{Hom}_{\mathfrak h'}(W^\vee,\indhg(V^\vee))
\stackrel{\sim}{\longrightarrow}
\operatorname{Diff}_{G'}( \mathcal V_X,\mathcal W_Y), 
\end{eqnarray}
or equivalently, 
\begin{equation*}
D_{X\to Y}:\operatorname{Hom}_{\mathfrak g'}(\operatorname{ind}_{\mathfrak{h}'}^{\mathfrak{g}'}(W^\vee),\indhg(V^\vee))
\stackrel{\sim}{\longrightarrow}
\operatorname{Diff}_{G'}( \mathcal V_X,\mathcal W_Y).
\leqno{(\ref{eqn:isoo})'}
\end{equation*}
\end{cor}

This result was previously known in the special setting where $X=Y$ and both are complex flag varieties,
 \emph{i.e.} $G' = G$ and $H' = H$ are Borel subgroups, see \emph{e.g.} \cite{Ks74, HJ}.

 When ${\mathfrak {g}}'$ is a reductive subalgebra
 and ${\mathfrak {h}}'$ is a parabolic subalgebra,
 the existence
 of an ${\mathfrak {h}}'$-module $W$
 for which the left-hand side of \eqref{eqn:isoo}
 is non-zero,
 is closely
 related to the \lq\lq{algebraic discrete decomposability}\rq\rq\
 of the ${\mathfrak {g}}$-module
 $\indhg (V^{\vee})$
 upon restriction to the subalgebra ${\mathfrak {g}}'$
 (\cite[ Part III]{xkdecomp}, \cite{K12}).  This relationship
plays a crucial role in the proof of the Localness Theorem (Theorem \ref{thm:C} below), which asserts that any continuous
 symmetry breaking operator in a holomorphic setting is necessarily realized as 
 a differential operator. \vskip8pt

An analogue of
Theorem \ref{thm:surject} holds in a holomorphic setting as well. 
More precisely,
 let $G_{\C}$ be a complex Lie group, $G_{\C}'$, $H_{\C}$
 and $H_{\C}'$ be closed complex subgroups
 such that $H_{\C}'\subset H_{\C} \cap G_{\C}'$. 
We write ${\mathfrak {g}}$, ${\mathfrak {h}}$, 
 $\ldots$
 for the Lie algebras
 of the complex Lie groups $G_{\C}$, $H_{\C}$, $\ldots$, respectively.
Given finite-dimensional holomorphic representations $V$ of $H_{\C}$
 and $W$ of $H_{\C}'$,
 we form the associated equivariant holomorphic vector bundles 
 \begin{eqnarray*}
\mathcal V:= G_{\C}\times_{H_\C} V&\mathrm{over}&
X_{\C}=G_{\C}/H_{\C},\\
\mathcal W:= G_{\C}'\times_{H_{\C}'}W&\mathrm{over}&Y_{\C}=G_{\C}'/H_{\C}'. 
\end{eqnarray*}

For simplicity,
 we assume
 that $H_{\C}'$ is connected.  
(This is always the case
 if $G_{\C}'$ is a connected complex reductive Lie group 
 and $H_{\C}'$ is a parabolic subgroup of $G_{\C}'$.)
Then we have:
\begin{thm}[Duality theorem in the holomorphic setting]\label{prop:2.10} 
There exists a canonical linear isomorphism:
$$
D_{X\to Y}:\operatorname{Hom}_{\mathfrak g'}\left(
\mathrm{ind}_{\mathfrak h'}^{\mathfrak g'}(W^\vee),
\mathrm{ind}_{\mathfrak h}^{\mathfrak g}(V^\vee)\right)
\stackrel{\sim}{\longrightarrow}
\operatorname{Diff}_{G_{\C}'}^{\mathrm{hol}}(\mathcal V_{X_{\C}},\mathcal W_{Y_{\C}}).
$$
\end{thm}

Suppose furthermore
 that $G$, $G'$, $H$ and $H'$ are real forms
 of the complex Lie groups
 $G_{\C}$, $G_{\C}'$, $H_{\C}$
 and $H_{\C}'$, 
 respectively.  
By restriction, we regard $V$ and $W$ as $H$- and $H'$-modules, and form the vector bundles 
 $$
 {\mathcal{V}}=G \times_H V\,\mathrm{over}\,\,X=G/H,\quad
 {\mathcal{W}}=G' \times_{H'} W\,\mathrm{over}\,\,Y=G'/H'.
 $$

One may then ask whether all symmetry breaking operators
 admit holomorphic extensions.  
Here is a simple sufficient condition:
\begin{cor}[{\cite[Corollary 2.13]{KP16a}}]
\label{cor:holoreal} 
If $H'$ is contained 
 in the connected complexification $H_{\C}'$, 
 then there exists a natural bijection:
\[
\operatorname{Diff}_{G_{\C}'}^{\operatorname{hol}}({\mathcal{V}}_{X_{\C}},{\mathcal{W}}_{Y_{\C}})
 \overset \sim \to 
\operatorname{Diff}_{G'}({\mathcal{V}}_{X},{\mathcal{W}}_{Y}).  
\]
\end{cor}

In general, symmetry breaking operators between two principal series representations of real reductive Lie groups $G'\subset G$
are given by integro-differential operators and their meromorphic continuations
in geometric models.  
Such operators were highlighted by Knapp and Stein in the case $G'=G$ and further developed by Kobayashi and Speh \cite{KS,KS18} in the broader setting $G'\subsetneqq G$. 
Among all symmetry breaking operators, those given by differential operators are quite rare, as evidenced by the complete classification results for certain pairs of real reductive groups $G'\subset G$.

However, when $X$ is a Hermitian symmetric space,
 $Y$ is a subsymmetric space, and $G'\subset G$ are the groups of biholomorphic transformations of $Y\hookrightarrow X$,
respectively, a remarkable ``inverse" phenomenon occurs. Namely, any continuous $G'$-intertwining operator between two representation spaces
consisting of holomorphic sections of holomorphic vector bundles over Hermitian symmetric spaces is necessarily realized as a differential operator. This is the second key result, which we call the \emph{Localness Theorem}.\vskip7pt

Let $G$ be a connected reductive Lie group, $\theta$ a Cartan involution, and $G/K$ the associated Riemannian symmetric space.
Let ${\mathfrak {c}}({\mathfrak {k}})$ 
denote the center
 of the complexified Lie algebra
 $\mathfrak k:=\operatorname{Lie}(K)\otimes_{\mathbb{R}} {\mathbb{C}}$.  
Suppose that $G/K$ is a Hermitian symmetric space. This means that
  there exists a characteristic element $Z \in 
{\mathfrak {c}}({\mathfrak {k}})$
 such that the eigenvalues
 of ${\operatorname{ad}} (Z) \in {\operatorname{End}}({\mathfrak {g}})$
 are 0 or $\pm 1$ and, that $\mathfrak {g}$ admits an eigenspace decomposition
\[
{\mathfrak {g}}={\mathfrak {k}}+{\mathfrak {n}}_+
 +{\mathfrak {n}}_-
\]
 with respect to ${\operatorname{ad}} (Z)$, corresponding to the
 eigenvalues $0$, $1$, and $-1$, respectively.  
We note that ${\mathfrak {c}}({\mathfrak {k}})$ is one-dimensional if $G$ is simple. 

 Let $G_\C$ be a complex reductive Lie group with Lie algebra $\mathfrak g$, and $P_\C$ the maximal parabolic subgroup with Lie algebra
 ${\mathfrak {p}}:= {\mathfrak {k}}+{\mathfrak {n}}_+$, 
 with abelian nilradical ${\mathfrak {n}}_+$. The complex structure
 of the homogeneous space $G/K$ is induced from
the open embedding
$$
G/K\subset G_\C/K_\C\exp\mathfrak n_+ = G_{\mathbb{C}}/P_{\mathbb{C}}.
$$

Let $G'$ be a connected reductive subgroup of $G$.
Without loss of generality we assume 
that $G'$ is $\theta$-stable. We set
 $K':=K\cap G'$.
 Furthermore, we assume that the characteristic element $Z\in\mathfrak c(\mathfrak k)$ satisfies
 \begin{equation}\label{eqn:ckk1}
 Z\in\mathfrak k'.
 \end{equation}
 
 If the condition \eqref{eqn:ckk1} holds, then the parabolic subalgebra $\mathfrak p$ is $\mathfrak g'$-compatible 
 %(that is  there exists a hyperbolic element $Z\in{\mathfrak {g}}^\prime$ such that ${\mathfrak {p}}={\mathfrak {p}}(Z)$), 
 and the homogeneous space
$G'/K'$ is a Hermitian sub-symmetric space of $G/K$ such that the embedding
$G'/K' \hookrightarrow G/K$ 
 is holomorphic.
 
Consider a finite-dimensional representation of $K$ on a complex vector space $V$. We extend it to a holomorphic representation of $P_\C$ 
 by letting the unipotent subgroup $\exp(\mathfrak n_+)$ act trivially,
 and form a holomorphic
 vector bundle $\mathcal V_{X_\C}=G_\C\times_{P_\C}V$ over $X_{\C}=G_\C/P_\C$.
 Restricting to the open subset $G/K\subset X_\C,$ we obtain a $G$-equivariant holomorphic vector bundle $\mathcal V:=G\times_K V$. We thus have
 a natural representation of $G$ on the vector space $\mathcal O(G/K,\mathcal V)$ of global holomorphic sections, endowed with
 the Fr\'echet topology of uniform convergence on compact sets.
 
 Likewise, given a finite-dimensional  representation $W$ of $K'$, we form
 the $G'$-equivariant holomorphic vector bundle $\mathcal W=G'\times_{K'}W$ over $G'/K'$ and consider
 the representation of $G'$ on $\mathcal O(G'/K',\mathcal W)$.
By definition, it is clear that
\begin{equation}\label{eqn:dcont}
\operatorname{Diff}_{G'}^{\mathrm{hol}}(\mathcal V_{X},\mathcal W_{Y}) 
\subset
\operatorname{Hom}_{G'}\left(\mathcal{O}\left(G/K,\mathcal V\right),\mathcal O(G'/K',\mathcal W)\right).
\end{equation}
The Localness Theorem  below shows that these two spaces in fact coincide.

\begin{thm}[Localness Theorem]\label{thm:C}
Let  $G'$ be a reductive
subgroup of $G$ satisfying \eqref{eqn:ckk1}, and
let $V$ and $W$ be any finite-dimensional  representations of $K$ and
$K'$, respectively. 
Then any continuous $G'$-homomor\-phism from
$\mathcal O(G/K,\mathcal V)$ to $\mathcal O( G'/K',\mathcal W)$ is
given by a holomorphic differential operator
 with respect to the holomorphic embedding between
  the Hermitian symmetric
spaces $G'/K'\hookrightarrow G/K$. Equivalently, 
$$\operatorname{Diff}_{G'}^{\mathrm{hol}}(\mathcal V_{X},\mathcal W_{Y}) 
=
\operatorname{Hom}_{G'}\left(\mathcal{O}\left(G/K,\mathcal V\right),\mathcal O(G'/K',\mathcal W)\right).
$$\end{thm}

This Localness Theorem reflects the algebraic discrete decomposability of the restriction of highest-weight modules of $\mathfrak g$ to the subalgebra $\mathfrak g'$.
Furthermore, this phenomenon extends to the complexification of underlying flag varieties, as formulated in
the following result, which constitutes the third key ingredient of our approach.

\begin{thm}[Extension Theorem]\label{thm:C'}
Under the assumptions of Theorem \ref{thm:C}, any differential operator in $\operatorname{Diff}_{G'}^{\mathrm{hol}}(\mathcal V_{X},\mathcal W_{Y})$
 (or, equivalently,
 any continuous $G'$-homomorphism)
 extends to a $G'_\C$-equivariant holomorphic differential operator
 with respect to a holomorphic map between the flag varieties 
 $$
 Y_{\C}=G'_\C/P_\C'\hookrightarrow  X_{\C}=G_\C/P_\C.
 $$
More precisely,
 the natural injection 
\begin{equation}\label{eqn:inj}
\operatorname{Diff}_{G'_\C}^{\mathrm{hol}}(\mathcal V_{X_{\C}},\mathcal W_{Y_{\C}}) 
\hookrightarrow
\operatorname{Diff}_{G'}^{\mathrm{hol}}\left(\mathcal V_{X},\mathcal W_{Y}\right)
\end{equation}
is bijective.
\end{thm}

It is noteworthy that the resulting holomorphic equivariant differential operators descend to other real forms of the pair $(G_\C, G_\C')$, thereby inducing a canonical injection of
$\operatorname{Diff}_{G'}^{\mathrm{hol}}\left(\mathcal V_{X},\mathcal W_{Y}\right)$
into the space of symmetry breaking operators on the corresponding real flag varieties. 
For example, the Rankin--Cohen brackets induce differential symmetry breaking operators between principal series representations of the groups $SL(2,\C)\supset SL(2,\R)$, see \cite{KoKuPe18} for example.
The converse statement is not true, since non-local symmetry breaking operators may occur.

\subsection{F-method}\label{sec:F-method}
Let $P'_\C=L'_\C\exp(\mathfrak n'_+)\subset G'_\C$ and $P_\C=L_\C\exp(\mathfrak n_+)\subset G_\C$ be compatible parabolic subgroups.
The $F$-method, encoded in the Duality Theorem \ref{thm:surject}, provides a powerful tool and an explicit procedure for determining differential symmetry breaking operator 
$$
D\in\mathrm{Diff}_{G'}(\mathcal V_X,\mathcal W_Y).
$$
More precisely, we put
$$
\varphi:=D_{X\to Y}^{-1}(D)\in\operatorname{Hom}_{H'}(W^\vee,\indhg(V^\vee)),
$$
and define its ``algebraic Fourier transform" $F_c(\varphi)$ as a $\operatorname{Hom}_\C(V,W)$-valued \emph{polynomial} on
$\mathfrak n_+$. This construction translates the equivariance condition satisfied by the operator $D$
into a (generally higher-order) system of partial differential equations satisfied by $F_c(\varphi)$. Consequently, finding a differential symmetry breaking operator $D$ is equivalent to finding a polynomial
$F_c(\varphi)$ satisfying a system of partial differential equations. This equivalence forms the essence of the F-method.

For simplicity, we consider the $F$-method in the particular setting where the nilradical $\mathfrak n'_+$ is abelian, that is, the symmetric space $G'/K'$ is Hermitian. 

In this situation, the algebraic Fourier transform $F_c(\varphi)$ agrees with the symbol of the corresponding differential operator $D_{X\to Y}(\varphi)$, and the method
can be summarized by the following diagram (see \cite[Sections 3 and 4]{KP16a}):
{\scriptsize
\begin{equation}\label{eqn:diagram}
\xymatrix{
&\left(\operatorname{Pol} (\mathfrak n_+)\otimes \operatorname{Hom}_\C(V,W)\right)^{L',\widehat{d\pi_\mu}(\mathfrak n_+')} & {}\\
\operatorname{Hom}_{P'}(W^\vee,\indpg(V^\vee)) \ar[ur]^{F_c}\ar[rr]_{D_{X\to Y}}^\sim
&{} & \operatorname{Diff}_{G'}(\mathcal V_X, \mathcal W_Y)\ar[ul]_{\mathrm{Symb}\otimes\mathrm{id}}
}
\end{equation}
}

\begin{remark}\label{rmk:n+}
The $F$-method remains valid even when $\mathfrak n_+'$ is not commutative. In this case, the diagram
\eqref{eqn:diagram} continues to apply, provided that the definition of the symbol map is appropriately modified.
\end{remark}

\subsection{Various examples}

Using the F-method, one can explicitly construct and classify differential symmetry breaking operators (differential SBOs)
in a variety of geometric settings $Y\subset X$. By the Duality Theorem \ref{thm:surject}, this problem is equivalent to an algebraic one, namely, the determination of \emph{singular vectors} of generalized Verma modules, as illustrated by an $\mathfrak{sl}_2$
example in Section \ref{sec:rec}.

Several recent examples include:
\begin{enumerate}
\item Rankin--Cohen brackets for $(X,Y)=(\mathbb H \times \mathbb H, \mathrm{diag}\,\mathbb H)$; see \cite{KP16}.
\item Juhl's conformally covariant SBOs  for $(X,Y)=(S^n,S^{n-1})$; see \cite{Juhl,KOSS}.
\item Extension of (2) to conformal geometry built on pseudo-Rie\-mannian manifolds $Y\subset X$; see \cite{KoKuPe18}.
\item Conformally covariant SBOs acting on \emph{differential forms}, rather than on functions, for $(X,Y)=(S^n,S^{n-1})$; see
\cite{Juhl,KoKuPe}.
\item Six geometries $(X,Y)$ corresponding to six pairs of Hermitian Lie groups $G\supset G'$ such that 
$\operatorname{rank}_\R G/G'=1$; see \cite{KP16}.
\end{enumerate}

Let us make a few comments on these examples.

In cases (1) and (2), both $V$ and $W$ are one-dimensional. Consequently, the symbols of the differential SBOs are given by scalar-valued polynomials, which are characterized via the F-method as polynomial solutions to certain second-order differential equations. By contrast, case (4) provides
the first instance of a classification of differential SBOs in which both $\dim V$ and $\dim W$ are greater than one.

Examples (2)--(4) serve as fundamental building blocks for the general classification theorem of SBOs between (degenerate) principal series representations when non-local SBOs are taken into account; see \cite{KoLe,KS,KS18}.

In the remainder of this section, we focus on cases (1) and (5) in detail.

In Section \ref{sec:blow}, we explain---through explicit computations based on the F-method---an intrinsic reason why the 
Rankin--Cohen brackets in case (1) are closely related to Jacobi polynomials, and somewhat surprisingly, why there exist
differential SBOs that \emph{are not of Rankin--Cohen type} when the parameters are highly singular. We analyze this phenomenon from
several complementary perspectives.

In Section \ref{sec:25} we turn to case (5), where we determine when a normal derivative defines an SBO, and subsequently formulate higher-dimensional analogues of Rankin--Cohen brackets.

\subsection{Blow-up of Multiplicities}\label{sec:blow}
In this section, we explain why and how the multiplicity-free property of fusion rules of Verma modules fails for highly singular values of spectral parameters, and we emphasize the role played by symmetry breaking operators with holomorphic parameters in this phenomenon.
In order to avoid unnecessary technicalities, we restrict our attention to the group $G=SU(1,1)$.

For $\lambda\in\Z$, we write $\mathcal L_\lambda$ for the $G$-equivariant holomorphic line bundle over the unit disk 
$$
D=\{z\in\C:\vert z\vert<1\}.
$$
Using the Bruhat decomposition, we trivialize the line bundle $\mathcal L_\lambda$ and identify
the regular representation of $G$ on
$\mathcal O(D,\mathcal L_\lambda)$ with the following multiplier representation on $\mathcal O(D)$, as in \eqref{eqn:hds} for
the action of $SL(2,\R)$ on $\mathcal H^2_\lambda(\mathbb H)$:
$$
\left(\pi_\lambda(g)F\right)(z)=(cz+d)^{-\lambda}F\left(\frac{az+b}{cz+d}\right),
\,\mathrm{for}\, g^{-1}=\begin{pmatrix}
a&b\\c&d
\end{pmatrix},
 F\in\mathcal O(D).
$$

For $\lambda\in\C$, the representation $\pi_\lambda$, originally defined for the discrete parameter
$\lambda\in\Z$, admits an extension to complex values of $\lambda$. In order to make sense of
this extension, it is necessary to pass from $G=SU(1,1)$ to the universal covering group $\widetilde G=\widetilde{SU}(1,1) $.

Let $\mathfrak b$ be a Borel subalgebra of $\mathfrak g=\mathfrak{sl}(2,\C)$ consisting of lower triangular matrices. The induced
module $\indbg(\nu)=U(\mathfrak g)\otimes_{U(\mathfrak b)}\C_\nu$
is a Verma module of the Lie algebra $\mathfrak g$.

In our parametrization, if $\lambda=1-k$ with $k\in\N$, then
the $k$-dimensional irreducible representation occurs both as a subrepresentation of
$(\pi_\lambda,\mathcal O(D))$ and as a quotient of $\indbg(-\lambda)$.

We consider symmetry breaking operators from the tensor product representation
$\mathcal O(\mathcal L_{\lambda'})\,\widehat\otimes\,\mathcal O(\mathcal L_{\lambda''})=
\mathcal O(D\times D,\mathcal L_{\lambda'}\boxtimes \mathcal L_{\lambda''})
$ to $\mathcal O(\mathcal L_{\lambda'''})$.
Here, $\widehat\otimes$ denote the completion of the tensor product of two nuclear spaces.
Theorem \ref{thm:U(n,1)} asserts that the Rankin--Cohen bracket $\mathcal{RC}_{\lambda',\lambda''}^{\lambda'''}$ provides an example of such an operator whenever 
$$
\lambda'''-\lambda'-\lambda''\in2\N.
$$ 
For generic values of the parameters, these operators exhaust all symmetry
breaking operators.
However, for certain highly singular values of spectral parameters $(\lambda',\lambda'',\lambda''')$, the Rankin--Cohen brackets vanish, whereas other differential symmetry breaking operators emerge.
In what follows, we examine this phenomenon in detail.

For $(\lambda',\lambda'',\lambda''')\in\C^3$, we define
$$
H(\lambda',\lambda'',\lambda'''):=\mathrm{Hom}_{\widetilde G}(\mathcal O(\mathcal L_{\lambda'})\widehat\otimes\mathcal O(\mathcal L_{\lambda''}), \mathcal O(\mathcal L_{\lambda'''})).
$$
By the Localness Theorem \ref{thm:C}, this space coincides with the space of differential symmetry breaking operators:
\begin{eqnarray*}
H(\lambda',\lambda'',\lambda''')&:=&
\mathrm{Diff}_{\widetilde G}(\mathcal O(\mathcal L_{\lambda'})\widehat\otimes\mathcal O(\mathcal L_{\lambda''}), \mathcal O(\mathcal L_{\lambda'''}))\\
&\simeq& \mathrm{Hom}_{\mathfrak g}(\indbg(-\lambda'''),\indbg(-\lambda')\otimes
\indbg(-\lambda'')),
\end{eqnarray*}
where the second isomorphism follows from the Duality Theorem \ref{thm:surject}. 

The general
theory of multiplicity-freeness of branching rules \cite{Kmf} implies that the dimension of $H(\lambda',\lambda'',\lambda''')$ is \emph{generically} equal to 0 or 1. However, finer properties of symmetry breaking operators allow us to control the behavior of multiplicities even beyond the generic case. 
The following result gives a precise dimension formula \cite[Theorem 9.1]{KP16}:

\begin{thm}\label{thm:Hlambdas}
The dimension of the vector space $H(\lambda',\lambda'',\lambda''')$ is uniformly bounded for all
$(\lambda',\lambda'',\lambda''')\in\C^3$. More precisely,
\begin{enumerate}
\item $\dim_\C H(\lambda',\lambda'',\lambda''')\in\{0,1,2\}$.
\item $H(\lambda',\lambda'',\lambda''')\neq\{0\}$ if and only if
\begin{equation}\label{eqn:leven}
\lambda'''-\lambda'-\lambda''\in 2\N.
\end{equation}
\item Suppose \eqref{eqn:leven} is satisfied. Then the following three conditions are equivalent:
\begin{enumerate}
\item[(i)] $\dim_\C H(\lambda',\lambda'',\lambda''')=2$.
\item[(ii)]
\begin{eqnarray}\label{eqn:l2}
&&\lambda',\lambda'',\lambda'''\in\Z,\quad 2\geq\lambda'+\lambda''+\lambda''',\quad\mathrm{and} \\
&&\lambda'''\geq\vert\lambda'-\lambda''\vert+2.\nonumber
\end{eqnarray}
\item[(iii)] $\mathcal{RC}_{\lambda',\lambda''}^{\lambda'''}=0$.
\end{enumerate}
\end{enumerate}
\end{thm}

A purely algebraic characterization of the multiplicity blow-up phenomenon has recently been formulated by R. Murakami \cite{mura24}, from the viewpoint of the extent to which complete reducibility fails for tensor product of Verma modules. More precisely, in the notations of Theorem \ref{thm:Hlambdas}, the condition $\mathcal{RC}_{\lambda',\lambda''}^{\lambda'''}=0$ is equivalent to
$$
(\mathrm{iv})\quad 2\leq \lambda'''\,\mathrm{and}\, \left(M(-\lambda')\otimes M(-\lambda'')\right)^{\chi_{-\lambda'''+1}}=
M(-\lambda''')\oplus M(\lambda'''+2).
$$
Here, $\chi_\nu: \mathcal Z(\mathfrak g)\To\C$ denotes the $\mathcal Z(\mathfrak g)$-infinitesimal
character via the Harish-Chandra isomorphism, and 
$M^{\chi_\nu}$ denotes the $\mathfrak{g}$-submodule of a module $M$ consisting of elements annihilated by a power of the ideal $\mathrm{Ker}(\chi_\nu:\mathcal Z(\mathfrak g)\to \C)$, that is the $\chi_\nu$-primary component of a locally $\mathcal Z(\mathfrak g)$-finite 
$\mathfrak g$-module $M$.

The original proof of Theorem \ref{thm:Hlambdas} relies on the 
F-method, which reduces problems in symmetry breaking to questions concerning differential equations. 

Conversely, one expects that the intriguing phenomena observed for symmetry breaking---such as
the blow-up of multiplicities---should also manifest themselves in the setting of differential equations.
Motivated by representation-theoretic considerations, we describe the corresponding phenomena for
polynomial solutions of the classical Jacobi differential equation and clarify their structure.

These results are derived purely from the classical theory of hypergeometric functions, including
Kummer's connection formul{\ae}, and we shall explain below how they relate to the structural theorems for symmetry breaking operators.
\vskip7pt

%In particular, the multiplicity-two phenomenon is reduced to the dimension formula for polynomial solutions of the Jacobi differential equation (Theorem \ref{thm:Hlambdas}, (2)), while the factorization of symmetry breaking operators stated in Theorem \ref{thm:Djs} below is reduced to the Kummer connection formulas, which are described in the following theorem.

Let $\ell\in\N$.  Define the finite set 
\begin{equation}\label{eqn:lmdl}
\Lambda_\ell:=\{(\alpha,\beta)\in\Z^2:\alpha+\ell\geq 0, \beta+\ell\geq0,\alpha+\beta\leq-(\ell+1)\},
\end{equation} which has the cardinality $\frac12\ell(\ell+1)$.
We can now formulate a precise characterization of the exceptional phenomenon for polynomial solutions of the Jacobi differential equation.

Here, the Jacobi differential equation refers to equation \eqref{eqn:JacobiDE}; throughout this section, the parameter $a$ appearing in \eqref{eqn:JacobiDE} is fixed to be $a=\ell$.

\begin{thm}[{\cite[Theorem 11.2]{KP16}}]\label{thm:Jacobizero}${}$
\begin{enumerate}
\item The following four conditions on $(\alpha,\beta)\in\C^2$ are equivalent:
\begin{enumerate}
\item[(i)] The Jacobi polynomial $P_\ell^{\alpha,\beta}(t)$ vanishes identically as a polynomial in $t$.
\item[(ii)] There exist two linearly independent polynomial solutions to the Jacobi differential equation
\eqref{eqn:JacobiDE}, of degree at most $\ell$.
 \item[(iii)] All solutions to the Jacobi differential equation are polynomial.
\item[(iv)] $(\alpha,\beta)\in\Lambda_\ell$. 
\end{enumerate}
\item If one (therefore any) of the equivalent conditions (i)-(iv) is satisfied, then
$\alpha$ is a negative integer and the limit
\begin{equation}\label{eqn:Fcnegative}
\lim_{\varepsilon\to0}{}_2F_1(-\ell,\alpha+\beta+1;\alpha+\varepsilon+1;z)
\end{equation}
exists and is a polynomial in $z$, which we denote by
$$
{}_2F_1(-\ell,\alpha+\beta+1;\alpha+1;z).
$$
Then any two of the following three polynomials
\begin{eqnarray}
g_1(z)&:=& z^{-\alpha} {}_2F_1(-\alpha-\ell, \beta+\ell+1; 1-\alpha;z),\label{eqn:g1}\\
g_2(z)&:=&{}_2F_1(-\ell,\alpha+\beta+\ell+1;\alpha+1;z),\label{eqn:g2}\\
g_3(z)&:=& (1-z)^{-\beta}{}_2F_1(-\beta-\ell,\alpha+\ell+1;1-\beta;1-z),\label{eqn:g3}
\end{eqnarray}
with $z=\frac12(1-t)$,
are linearly independent polynomial solutions to the Jacobi differential equation \eqref{eqn:JacobiDE}
of degree 
$$
\ell, \quad -(\alpha+\beta+\ell+1), \quad {\mathrm{and}}\quad \ell,
$$ 
respectively. In particular, any polynomial solution to \eqref{eqn:JacobiDE} is of degree at most $\ell$.
\end{enumerate}
\end{thm}

We illustrate the phenomena occurring for the exceptional parameter space $\Lambda_\ell$ from 
two complementary viewpoints: the Jacobi differential equation and representation theory.
The correspondence is summarized in the following table:

\begin{table}[!h]
\begin{center}
\begin{tabular}{c|c}
Differential Equations & Representation Theory\\
\hline
{}&{}\\
Jacobi Polynomial vanish & Rankin--Cohen brackets vanish\\
{\footnotesize(Theorem \ref{thm:Jacobizero} (1))} & {\footnotesize(Theorem \ref{thm:Hlambdas}(3))}\\
\hline
All solutions are polynomial & Blow-up of multiplicities\\
{\footnotesize(Theorem \ref{thm:Jacobizero} (1))} & {\footnotesize(Theorem \ref{thm:Hlambdas}(2))}\\\hline
Kummer's connection formul{\ae}  & Factorization of SBOs\\
with vanishing coefficients&\\
{\footnotesize(Theorem \ref{thm:Jacobizero} (2))} & {\footnotesize(Theorem \ref{thm:Djs})}
\end{tabular}
\end{center}
\label{default}
\end{table}%

To elucidate the relationship with the \emph{factorization identities} of symmetry breaking operators and the special solutions of the Jacobi differential equation appearing in the bottom row of the table, we express
the polynomials $g_j(z)$ ($j=1,2,3$) in Theorem  \ref{thm:Jacobizero} in terms of the Jacobi polynomials.
Suppose $(\alpha,\beta)\in\Lambda_\ell$. Then, the three polynomials $g_1,g_2,g_3$ introduced in Theorem \ref{thm:Jacobizero} admit the following renormalizations, which can be expressed in terms of Jacobi polynomials with various parameters:
\begin{itemize}
\item[(1)] Define the renormalized polynomials by
\begin{align*}
\widetilde g_1(z)&:=
\binom{\ell}{-\alpha}\cdot g_1(z)&&=z^{-\alpha}P_{\ell+\alpha}^{-\alpha,\beta}(1-2z);\\
\widetilde g_2(z)&:=(-1)^{-\ell-\alpha-\beta-1}\binom{-\alpha-1}{\ell+\beta}\cdot g_2(z)&&=P_{-\ell-\alpha-\beta-1}^{\alpha,\beta}(1-2z);\\
\widetilde g_3(z)&:= (-1)^{\beta+\ell}\binom{\ell}{-\beta}\cdot g_3(z)&&=(1-z)^{-\beta}P_{\ell+\beta}^{\alpha,-\beta}(1-2z).
\end{align*}
\item[(2)] Their linear dependence is explicitly given by 
$$
(-1)^\alpha\,\widetilde g_3(z)=\widetilde g_1(z)-\widetilde g_2(z),
$$
that is 
$$
P_{-\ell-\alpha-\beta-1}^{\alpha,\beta}(t)=
(-1)^{\alpha+1}\left(\frac{1+t}2\right)^{-\beta}P_{\beta+\ell}^{\alpha,-\beta}(t)
+\left(\frac{1-t}2\right)^{-\alpha}P_{\alpha+\ell}^{-\alpha,\beta}(t).
$$
\end{itemize}

By means of this identity, we obtain an explicit basis of $H(\lambda',\lambda'',\lambda''')$
with a transparent interpretation from the viewpoint of representation theory. Let
\begin{equation}\label{eqn:abclmd}
\alpha:=\lambda'-1,\quad \beta:=1-\lambda''',\quad \ell:=\frac12(\lambda'''-\lambda'-\lambda'').
\end{equation}
As in \eqref{eqn:Jacobi_two}, we inflate the polynomials $\widetilde g_j$ to homogeneous polynomials of degree $\ell$ in two variables by
$$
G_j(x,y):=(-y)^\ell \widetilde g_j\left(1+\frac{2x}y\right),
$$ 
and define bi-differential operators
$$
D_j:=\mathrm{Rest}_{z_1=z_2=z}\circ G_j\left(\frac{\partial}{\partial z_1},\frac{\partial}{\partial z_2}\right),
\quad j=1,2,3.
$$

The following factorization result for symmetry breaking operators clarifies what happens with the triple of representations $\pi_{\lambda'},\pi_{\lambda''},\pi_{\lambda'''}$ when the blow-up
phenomenon of multiplicities occurs. This  is the representation-theoretic counterpart of the vanishing of coefficients in Kummer's connection formula for these specific parameters.

\begin{thm}\label{thm:Djs}
Suppose that the conditions \eqref{eqn:leven} and \eqref{eqn:l2} are satisfied. 
Then :
\begin{enumerate}
\item The operators $D_j$ ($j=1,2,3$) are non-zero $G$-homomorphisms from $\mathcal O(\mathcal L_{\lambda'})\widehat\otimes \mathcal O(\mathcal L_{\lambda''})$ to $\mathcal O(\mathcal L_{\lambda'''})$.
\item
None of the operators $D_j$ ($j=1,2,3$) is a Rankin--Cohen operator. Furthermore,
 $$
 1-\lambda', 1-\lambda''\quad {\mathrm {and}}\quad
\lambda'''-1\in\N_+,
$$ 
and each $D_j$ ($j=1,2,3$) factorizes into a composition of two natural intertwining operators as follows:
\begin{eqnarray*}
D_1&=&\mathcal{RC}_{2-\lambda',\lambda''}^{\lambda'''}\circ
\left(\left(\frac{\partial}{\partial z_1}\right)^{1-\lambda'}\otimes\mathrm{id}\right),\\
D_2&=&
\mathcal{RC}_{\lambda',2-\lambda''}^{\lambda'''}\circ\left(
\mathrm{id}\otimes\left(\frac{\partial}{\partial z_2}\right)^{1-\lambda''}\right),\\
D_3&=&
\left(\frac{d}{dz}\right)^{\lambda'''-1}\circ
\mathcal{RC}_{\lambda',\lambda''}^{2-\lambda'''}.
\end{eqnarray*}
\item
The following linear relation holds:
\begin{eqnarray*}
D_1-D_2+(-1)^{\lambda'}D_3=0.
\end{eqnarray*}
\end{enumerate}
\end{thm}
The factorizations  in Theorem \ref{thm:Djs} are illustrated by the following diagram:
{\scriptsize
\begin{equation}\label{eqn:diagram3}
\xymatrix{
 & &\qquad\quad\mathcal O(\mathcal L_{2-\lambda'})\, \widehat\otimes\,
\mathcal O(\mathcal L_{\lambda''})\ar[drr]^{\mathcal{RC}_{2-\lambda',\lambda''}^{\lambda'''}}\qquad\qquad&  &\\
\mathcal O(\mathcal L_{\lambda'}) \widehat\otimes
\mathcal O(\mathcal L_{\lambda''})\ar[urr]^{\left(\frac{\partial}{\partial z_1}\right)^{1-\lambda'}\otimes\,\mathrm{id}}
\ar[rr]^{\mathrm{id}\,\otimes\left(\frac{\partial}{\partial z_2}\right)^{1-\lambda''}}
\ar[drr]_{\mathcal{RC}_{\lambda',\lambda''}^{2-\lambda'''}} 
& & \mathcal O(\mathcal L_{\lambda'}) \,\widehat\otimes\,
\mathcal O(\mathcal L_{2-\lambda''})\ar[rr]^{\mathcal{RC}_{\lambda',2-\lambda''}^{\lambda'''}} & &\mathcal O(\mathcal L_{\lambda'''}),\\
& &\mathcal O(\mathcal L_{2-\lambda'''})\ar[urr]_{\left(\frac{d}{d z}\right)^{\lambda'''-1}}& &
}
\end{equation}
}
Including these special case, we provide explicit bases of the space
$H(\lambda',\lambda'',\lambda''')$ of symmetry breaking operators, by dividing the parameters setting into
the following three mutually exclusive cases.\vskip7pt
\begin{enumerate}
\item[Case 0.] $\lambda'''-\lambda'-\lambda''\not\in 2\N$.
\item[Case 1.] $\lambda'''-\lambda'-\lambda''\in 2\N,$ but the condition \eqref{eqn:l2} is not fulfilled.
\item[Case 2.] $\lambda'''-\lambda'-\lambda''\in 2\N,$ and the condition \eqref{eqn:l2} is satisfied.
\end{enumerate}

\begin{cor}\label{cor:83}
$$
H(\lambda',\lambda'',\lambda''')=\left\{
\begin{matrix}
\{0\} & \mathrm{Case}\,0,\\
\C\cdot\mathcal{RC}_{\lambda',\lambda''}^{\lambda'''}&\mathrm{Case}\, 1,\\
\C D_1\oplus\C D_2=\C D_1\oplus\C D_3=\C  D_2\oplus\C D_3&\mathrm{Case}\, 2.
\end{matrix}\right.
$$
\end{cor}

\subsection{Higher-dimensional examples}\label{sec:25}
Explicit formul\ae\, for all differential symmetry breaking
operators were obtained by applying the above method for all six complex parabolic geometries $X\supset Y$ arising from
pairs of Hermitian Lie groups $(G,G')$ such that the semisimple symmetric space $G/G'$ is of split rank one (even though $G/G'$ itself may have higher rank); see \cite{KP16}.  

This approach also provides an intrinsic reason for why the coefficients of certain families of orthogonal polynomials appear as coefficients in explicit expressions of symmetry breaking operators of Rankin--Cohen type
(e.g., Rankin--Cohen brackets)
in three of the geometries, while normal derivatives define symmetry breaking operators in the remaining three cases. 

More precisely, the latter geometries $Y\hookrightarrow X$ correspond to the following inclusions of flag varieties:
\begin{eqnarray*}
\mathrm {Gr}_{p-1}(\C^{p+q})&\hookrightarrow& \mathrm {Gr}_{p}(\C^{p+q}),
\\
\mathbb P^n\C&\hookrightarrow&\mathrm{Q}^{2n}\C,
\\
\mathrm{IGr}_{n-1}(\C^{2n-2})&\hookrightarrow & \mathrm{IGr}_{n}(\C^{2n}),
\end{eqnarray*}
where $\mathrm {Gr}_{p}(\C^{k})$ denotes the complex Grassmanian of $p$-planes in $\C^k$,
 $\mathrm Q^{m}\C:=\{z\in\mathbb P^{m+1}\C\,:\, z_0^2+\cdots+z_{m+1}^2=0\}$ is the complex quadric,
 and $\mathrm{IGr}_{n}(\C^{2n}):=\{V\subset\C^{2n}\,:\,\mathrm{dim}\,V=n,
 \,Q|_V\equiv0\}$ is the Grassmanian of isotropic subspaces of $\C^{2n}$ with respect to a non-degenerate quadratic form $Q$.
 The corresponding three pairs of complexified reductive groups $(G_\C,G'_\C)$ are
\begin{eqnarray*}
(GL(p+q,\C), GL(p+q-1,\C)\times GL(1,\C)),\\
(SO(2n+2,\C), GL(n+1,\C))\\
(SO(2n+2,\C), SO(2n,\C)\times SO(2,\C)).
\end{eqnarray*}
On the other hand, differential symmetry breaking operators arising from the other three geometries include higher-dimensional analogues of Rankin--Cohen brackets, which we now explain.

From now on, as a representative case, we consider the pair 
$$
(G_\C,G'_\C)=(GL(n+1,\C)\times GL(n+1,\C), GL(n+1,\C)),
$$
 where the real form $G'$ is taken to be the Hermitian Lie group $U(n,1)$. 
The case $n=1$ was essentially treated in the previous section.
Our method provides a conceptual understanding of higher-dimensional generalization of Rankin--Cohen brackets as follows.

Consider the unit ball in $\C^n$, defined by 
$$
D:=\{Z\in\C^n\,:\, \Vert Z\Vert<1\},
$$
 where 
 $\Vert Z \Vert^2:=\sum_{j=1}^n |z_j|^2$
 for $Z=(z_1, \cdots, z_n)$.
This is the Hermitian symmetric domain of type AIII in $\C^n$ in the \'E. Cartan classification. 
As in the $n=1$ case treated in Section \ref{sec:blow}, the Lie group $G'=U(n,1)$ acts biholomorphically on $D$ by linear fractional transformations:
$$
g\cdot Z=(aZ+b)(cZ+d)^{-1}\quad{\mathrm{for}}\quad
g=\begin{pmatrix}a&b\\c&d\end{pmatrix}\in U(n,1),\,Z\in D.
$$
The isotropy subgroup at the origin is isomorphic to 
$K'=U(n)\times U(1)$. 

For $\lambda_1,\lambda_2\in\Z$, we form 
a $U(n,1)$-equivariant holomorphic line bundle $\displaystyle\mathcal L_{(\lambda_1,\lambda_2)}=U(n,1)\times_{U(n)\times U(1)}\C_{(\lambda_1,\lambda_2)}$ over $D$ associated with the holomorphic character of the complexified isotropy group $K'_\C=GL(n,\C)\times GL(1,\C)$:
$$
GL(n,\C)\times GL(1,\C)\To \C^\times, \quad (A,d)\mapsto (\det A)^{-\lambda_1}d^{-\lambda_2}.
$$
 
The representation of $U(n,1)$ on $\mathcal O(D,\mathcal L_{(\lambda_1,\lambda_2)})$ can be identified, via the Bruhat open cell trivialization, with the multiplier representation, denoted simply by $\pi_{(\lambda_1,\lambda_2)}$, of $U(n,1)$ on $\mathcal O(D)$ given by
$$
\left(\pi_{(\lambda_1,\lambda_2)}(g)F\right)(Z)=(cZ+d)^{-\lambda_1+\lambda_2}(\det g)^{-\lambda_1} F\left(({aZ+b})({cZ+d})^{-1}\right).
$$

Next, recall from \eqref{eqn:Jacobi_two} the definition of the  homogeneous polynomial in two variables obtained by inflating the classical Jacobi polynomial: 
\begin{equation}\label{eqn:inflJacobi}
P_\ell^{\alpha,\beta}(x,y):=y^\ell P_\ell^{\alpha,\beta}\left(2\frac xy+1\right).
\end{equation}
For instance, $P_0^{\alpha,\beta}(x,y)=1$, $P_1^{\alpha,\beta}(x,y)=(2+\alpha+\beta)x+(\alpha+1)y$, etc.
\vskip7pt

Let $\widetilde{U}(n,1)$ denote the universal covering of the group $U(n,1)$. Then, for
all $\lambda_1,\lambda_2\in\C$ one defines a 
$\widetilde{U}(n,1)$-equivariant holomorphic line bundle $\mathcal L_{(\lambda_1,\lambda_2)}$ over $D$, and the corresponding representation of $\widetilde{U}(n,1)$ on $\mathcal O(D,\mathcal L_{(\lambda_1,\lambda_2)})$.

In contrast to the $n=1$ case, representations on the spaces of holomorphic sections of line bundles are no longer closed under tensor product decompositions; that is, most of the irreducible components that appear in the decomposition cannot be realized as spaces of sections of line bundles (see \cite{Kmf} for explicit branching laws in the general setting). This naturally leads one to consider holomorphic vector bundles, rather than line bundles alone, as the target of symmetry breaking operators.

For this purpose, let $a\in\N, \lambda_1,\lambda_2\in\Z$, and let $W_{(\lambda_1,\lambda_2)}^a$
denote the irreducible representation of $K'_\C=GL(n,\C)\times GL(1,\C)$ having the highest weight
$(-\lambda_1,\cdots,-\lambda_1,-\lambda_1-a;-\lambda_2+a)$. Then $W_{(\lambda_1,\lambda_2)}^a$
can be realized on the space $\operatorname{Pol}^a[v_1,\cdots,v_n]$ of homogeneous polynomials of degree $a$, where $(v_1,\cdots,v_n)$ are the standard coordinates on
$(\C^n\times\C^n)/\mathrm{diag}\C^n\simeq\C^n$.

Differential operators on $D\times D$ with values in the vector space $W_{(\lambda_1,\lambda_2)}^a$
can be written as elements in
$$
\C\left[\frac{\partial}{\partial z'_1},\cdots,\frac{\partial}{\partial z'_n},\frac{\partial}{\partial z''_1},\cdots
\frac{\partial}{\partial z''_n}\right]\otimes\operatorname{Pol}^a[v_1,\cdots,v_n].
$$
Associated with the $K_\C'$-module $W_{(\lambda_1,\lambda_2)}^a$, we form a $G'$-equivariant holomorphic vector bundle $\mathcal W_{(\lambda_1,\lambda_2)}^a$ over $D=G'/K'$.

As in the line bundle case, for all $\lambda_1,\lambda_2\in\C$, one can define a corresponding
$\widetilde U(n,1)$-equivariant vector bundle, which we also denote $\mathcal W_{(\lambda_1,\lambda_2)}^a$. Similarly to $\pi_{(\lambda_1,\lambda_2)}$ for the line
bundle $\mathcal L_{(\lambda_1,\lambda_2)}$, we define the representation 
$\pi_{W_{(\lambda_1,\lambda_2)}^a}$ of $\widetilde {G'}$ on $\mathcal O(D,\mathcal W_{(\lambda_1,\lambda_2)}^a)
\simeq\mathcal O(D)\otimes \mathcal W_{(\lambda_1,\lambda_2)}^a$.

We are now ready to describe symmetry breaking operators in the higher-dimensional case.

\begin{thm}[\cite{KP16}, Theorem 8.1]\label{thm:U(n,1)}
Let $\lambda_1',\lambda_2', \lambda_1'',\lambda_2''\in\C, a\in \N$, and
 $\lambda':=\lambda_1'-\lambda_2'$
 and $\lambda'':=\lambda_1''-\lambda_2''$. 
\begin{enumerate}
\item The space
 $$
 \mathrm{Hom}_{\widetilde{U}(n,1)}\left(
 {\mathcal{O}}
(D, {\mathcal{L}}_{(\lambda_1',\lambda_2')} )\,\widehat\otimes \,\mathcal{O}
(D, {\mathcal{L}}_{(\lambda_1'',\lambda_2'')}),{\mathcal{O}}(D, {\mathcal{W}}_{(\lambda_1'+\lambda_1'', \lambda_2'+\lambda_2'')}^a)\right)
$$
has dimension either one or two. It is equal to two if and only if
\begin{equation}\label{eqn:Udim}
\lambda',\lambda''\in\{-1,-2,\cdots\}\quad\mathrm{and}\quad a\geq \lambda'+\lambda''+2a-1\geq\vert\lambda'-\lambda''\vert.
\end{equation}
\item The vector-valued differential
operator %from $\mathcal O(D\times D)$ to $\mathcal O(D)\otimes
%\mathrm{Pol}^a[v_1,\cdots,v_n]$ defined by
\begin{equation}\label{eqn:DXYUN1}
 D_{X\to Y,a}:=
  P_a^{\lambda'-1,-\lambda'-\lambda''-2a+1}\left(
\sum_{i=1}^{n}v_i \frac{\partial}{\partial z_i},
\sum_{j=1}^{n}v_j \frac{\partial}{\partial z_j}
\right)
\end{equation}
intertwines $\pi_{(\lambda_1',\lambda_2')}\boxtimes \pi_{(\lambda_1'',\lambda_2'')}\Big\vert_{G'}$ with $\pi_W$, where $W\simeq W^a_{(\lambda_1'+\lambda_1'',\lambda_2'+\lambda_2'')}$.

\item If the parameters $(\lambda',\lambda'',a)$ satisfy the condition \eqref{eqn:Udim}, then $ D_{X\to Y,a}=0$. Otherwise,
any symmetry breaking operator  
is proportional to $D_{X\to Y,a}$.
\end{enumerate}
\end{thm}

\begin{rem}
${}$
\begin{enumerate}

\item Part (1) of Theorem \ref{thm:U(n,1)} shows that, as in the $n=1$ case, the blow-up phenomenon of multiplicities can also occur in higher dimensions. The condition \eqref{eqn:Udim} specifies precisely when this occurs. As noted in the previous section, this phenomenon exhibits a rich interplay not only within representation-theoretic framework but also with other areas of mathematics, such as special solutions of differential equations.

\item The vector bundle
 ${\mathcal{W}}_{(\lambda_1, \lambda_2)}^{a}$
 is a line bundle
 if and only if $a=0$
 or $n=1$.  
In the case $n=1$, 
the formula \eqref{eqn:DXYUN1} reduces to
the classical Rankin--Cohen bidifferential operators (see \eqref{eqn:RCB}) with an appropriate
choice of spectral parameters, namely, for $a:=\frac12(\lambda'''-\lambda'-\lambda'')\in\N$, the following identity holds:
 \begin{equation}\label{eqn:JacobiRankin}
 \mathcal{RC}_{\lambda',\lambda''}^{\lambda'''}=(-1)^a
P^{\lambda'-1, 1-\lambda'''}_a
\left(\frac{\partial}{\partial z_1}, \frac{\partial}{\partial z_2}\right)\big\vert_{z_1=z_2=z}.
 \end{equation}
 
 \end{enumerate}
 \end{rem}
 \begin{rem}
\begin{enumerate}
Let $X=D\times D\supset Y=D$.
 \item If  $\lambda_1',\lambda_2', \lambda_1'',\lambda_2''\in\Z$ and $a\in\N$, then the 
 linear groups $G$ and $G'$ act equivariantly on the two bundles $\mathcal L_{(\lambda_1',\lambda_2')}\boxtimes 
\mathcal L_{(\lambda_1'',\lambda_2'')}\to D\times D$ and $\mathcal W_{(\lambda_1,\lambda_2)}^a\to D$, respectively.
 
\item If $\lambda',\lambda''>n$, then analogous statements as in Theorem \ref{thm:U(n,1)} remain true
 for continuous $G'$-homomorphisms
 between the Hilbert spaces
${\mathcal{H}}^2
\left(X, {\mathcal{L}}_{(\lambda_1',\lambda_2')} \boxtimes {\mathcal{L}}_{(\lambda_1'',\lambda_2'')})\right)$
 and 
${\mathcal{H}}^2
\left(Y, {\mathcal{W}}^a_{(\lambda_1'+\lambda_1'' , \lambda_2'+\lambda_2'')}\right)$.
\item Similar statements hold for continuous $G'$-homomorphisms
  between the Casselman--Wallach globalizations
 by the Localness Theorem
 \ref{thm:C}.  
\end{enumerate}
\end{rem}

\vskip12pt

\section{Further Perspectives}
Building on a detailed understanding of a broad range of symmetry breaking operators, we initiated in \cite{KP} a new line of investigation
 into branching problems
 by introducing the concept of {\emph{generating operator}}. This operator is defined
 as a formal power series in $t$,
  \begin{equation}
\label{eqn:Rlt}
T(t)\equiv T(\{R_{\ell}\}; t):=\sum_{\ell=0}^{\infty}\frac{R_{\ell}}{\ell!}t^{\ell} \in \operatorname{Hom}(\Gamma(X), \Gamma(Y)) \otimes {\mathbb{C}}[[t]],
\end{equation}
 for a family of operators $R_{\ell}: \Gamma(X)\To \Gamma(Y)$, where $X$ and $Y$ are smooth manifolds, and $\Gamma(X)$ and $\Gamma(Y)$ denote appropriate function (or section) spaces defined on them. 
 
 In the case when $X=\{\text{point}\}$, each
 $R_{\ell}$ may be identified with an element of $\Gamma(Y)$, 
 and the generating operator reduces to the classical notion of a {\it{generating function}}.
 Such functions have played a prominent role
 in the classical theory of orthogonal polynomials, especially
when $\Gamma(Y)={\mathbb{C}}[y]$.  

When $X=Y$, the space
 $\operatorname{Hom}(\Gamma(X), \Gamma(Y))\simeq {\operatorname{End}}(\Gamma(X))$
 carries a natural ring structure,
 and one may take $R_{\ell}$
 to be the $\ell$-th power of a single operator $R$
 on $X$.  
In this case, 
 the operator $T(t)$ in \eqref{eqn:Rlt}
 may be written as $e^{t R}$,
 provided that the series converges.  
We note that even if $R$ is a differential operator 
 on a manifold $X$, 
 the resulting operator $T(t)=e^{t R}$ is generally no longer differential.  
For example, 
 if $R=\frac {d}{d z}$ acts on ${\mathcal{O}}({\mathbb{C}})$, 
 then $T(t)=e^{t \frac d {dz}}$ is the shift operator
 $f(z) \mapsto f(z+t)$.  

In the very special setting where $X=Y$ and $R$ is a self-adjoint operator 
 with spectrum bounded from above, 
 the generating operator $T$ reduces to the classical
  {\it{semigroup}} $e^{t R}$, which is
 generated by $R$ for $\operatorname{Re} t>0$ and has been extensively studied.
 Typical examples include 
 the heat kernel
 for $R=\Delta$, 
the Hermite semigroup
 for $R=\frac 1 4 (\Delta-|x|^2)$  on $L^2({\mathbb{R}}^n)$, 
or
the Laguerre semigroup
 for $R=|x|(\frac{\Delta}4-1)$ on $L^2({\mathbb{R}}^n, \frac{1}{|x|} d x)$.  

We propose to consider a more general setting
 where we allow $X \ne \{\text{point}\}$ and $X \ne Y$.  
In this broader context, 
 we refer to $T$ in \eqref{eqn:Rlt}
 as the {\it{generating operator}}
 for the family of operators
 $R_{\ell} \colon \Gamma(X) \to \Gamma(Y)$.

  More specifically, when $(X,Y)= (D\times D, \mathrm{diag}\,D)$, we
  found, after choosing an appropriate normalization, an explicit closed formula for the generating operator
 associated with the family of Rankin--Cohen operators $R_{\ell}=\mathcal{RC}_{\lambda',\lambda''}^{\lambda'+\lambda''+2\ell}$   (\cite{KPshort, KP}).  This concept aims to control branching rules in an integrated manner through a global analytic object.\vskip7pt

As an application of the concept of generating operator, 
 we proposed in \cite{K23} a method 
 for transferring {\emph{discrete data}}
 into {\emph{continuous data}} via the generating operator. 
 This is accomplished by solving a Mittag-Leffler type problem for operator-valued meromorphic functions. We
 illustrated this method with an example involving $S L_2$ demonstrating that
the {\it{generating operator}} $T$ 
associated with the Rankin--Cohen brackets $\{R_{\ell}\}_{\ell \in {\mathbb{N}}}$
yields various families
of non-local intertwining operators with continuous parameters, which emerge
in geometric settings distinct from the original domains of the operators.  
 
Specifically, the countable family of holomorphic differential operators 
 $\{R_{\ell}\}_{\ell \in {\mathbb{N}}}$ initially defined
 on the Poincar\'e upper half-plane, leads to a generating operator for which the formal power series converges near $t=0$.
Although the generating operator itself is non-local, its generator can be extended to the entire complex projective space via the Extension Theorem \ref{thm:C'}. Consequently, the generating operator is also well-defined on the complex projective space. To recover the original discrete data (the Rankin--Cohen brackets), one differentiate with respect to $t$. Alternatively, by applying
a \emph{fractional derivative} in $t$, one can obtain continuous data---specifically, symmetry breaking operators with continuous parameters.

 When restricted to a different real form, it produces a family of covariant operators 
 in various geometric settings, including symmetry breaking for tensor products of principal series representations of $SL(2,\R)$, Poisson transforms for
  the de Sitter space $d\mathbb S^2$, and concrete embedding of discrete series representations into principal series representations. 
  
  This strategy is also relevant for higher-rank groups: extensions of generating operators associated with families of normal derivatives on Hermitian symmetric spaces give rise to new differential symmetry breaking operators on real Grassmannians.
\vskip7pt

We believe that this represents a challenging and promising research direction.
We expect it to enhance our understanding of 
branching problems for infinite-dimensional representations of higher-rank reductive Lie groups.

\section*{Acknowledgments}
The authors express their deepest gratitude to the organizers of the Thematic Program
\emph{Representation Theory and Noncommutative Geometry}, Alexandre Afgoustidis,
Anne-Marie Aubert,
Pierre Clare,
Jan Frahm,
Angela Pasquale,
Haluk \c{S}eng\"un  for the invitation and
acknowledge support of the Institut Henri Poincar\'e (UAR 839 CNRS-Sorbonne Universit\'e), and LabEx CARMIN (ANR-10-LABX-59-01)
to this Program.
\vskip7pt

The first author was partially supported by the JSPS under the Grant-in Aid for Scientific Research
(A) (JP23H00084).

\vskip 3pc

%\leftline
%{\sc {Michael Pevzner}}
%\leftline
%{
%Universit{\'e} de Reims-Champagne-Ardenne, 
%CNRS UMR 9008, LMR}
%\leftline{Reims, France} 
%
%\&
%
%\leftline
%{French-Japanese Laboratory
% in Mathematics and its Interactions, }
%\leftline{
%FJ-LMI CNRS IRL 2025, Tokyo, Japan}
%\leftline{E-mail: \texttt{ pevzner@math.cnrs.fr}}

\end{document}